\documentclass[11pt,a4paper,reqno]{amsart}
\usepackage[T1]{fontenc}
\usepackage{amsmath}
\usepackage{amsthm}
\usepackage{amsfonts}
\usepackage{amssymb}
\usepackage{graphicx}
\usepackage{amsbsy}
\usepackage{mathrsfs}
\usepackage{bbm}
\newcommand{\be}{\begin{equation}}
\newcommand{\ee}{\end{equation}}
\newcommand{\ba}{\begin{array}}
\newcommand{\ea}{\end{array}}
\newcommand{\bea}{\begin{eqnarray}}
\newcommand{\eea}{\end{eqnarray}}
\newcommand{\bee}{\begin{eqnarray*}}
\newcommand{\eee}{\end{eqnarray*}}

\newtheorem{Thm}{Theorem}[section]
\newtheorem{Lemma}[Thm]{Lemma}
\newtheorem{Prop}[Thm]{Proposition}

\newtheorem{Cor}[Thm]{Corollary}
\newtheorem{remark}[Thm]{Remark}
\numberwithin{equation}{section}

\catcode`@=11
\def\section{\@startsection{section}{1}%
  \z@{1.5\linespacing\@plus\linespacing}{.5\linespacing}%
  {\normalfont\bfseries\large\centering}}
\catcode`@=12

\def\C{{\mathbb C}}
\def\RR{\mathbb{R}}
\def\R{{\mathbb R}}

\def\lim{\mathop{\rm lim}}
\def\goto{\rightarrow}

\def\sup{\mathop{\rm sup}}

\def\e{\varepsilon}
\def\l{\lambda}

\def\log{{\rm log}}

\def\ue{\underline{\e}}
\def\tu{\tilde{u}}

\def\lsl{\frac{\lambda_s}{\lambda}}

\def\tgamma{{\tilde{\gamma}}}

\def\p{\partial}

\def\n{\nabla}
\def\l{\lambda}
\def\half{\frac 12}

\def\S{\Sigma}
\def\T{\Theta}

\def\ut{\tilde{u}}

\def\fref#1{{\rm (\ref{#1})}}

\def\ds{\displaystyle}

\def\pa{\partial}

\def\qte{\tilde{Q}_{\eta}}

\def\tq{\tilde{Q}}
\def\a{\alpha}

\def\uv{\underline{v}}
\def\ttu{\tilde{\tilde{u}}}
\def\ts{\tilde{\Sigma}}
\def\tt{\tilde{\Theta}}

\title[]{The instability of Bourgain-Wang solutions for the $L^2$ critical NLS}
\author[F. Merle]{Frank Merle}
\address{Universit\'e de Cergy Pontoise and IHES, France}
\email{frank.merle@math.u-cergy.fr}
\author[P. Rapha\"el]{Pierre Rapha\"el}
\address{IMT, Universit\'e Paul Sabatier, Toulouse, France}
\email{pierre.raphael@math.univ-toulouse.fr}
\author[J. Szeftel]{Jeremie Szeftel}
\address{DMA, Ecole Normale Sup\'erieure, France}
\email{szeftel@dma.ens.fr}

\begin{document}

\begin{abstract} We consider the two dimensional $L^2$ critical nonlinear Schr\"odinger equation  $i\pa_tu+\Delta u+u|u|^2=0$. In the pioneering work \cite{BW}, Bourgain and Wang have constructed smooth solutions which blow up in finite time $T<+\infty$ with the pseudo conformal speed $$\|\nabla u(t)\|_{L^2}\sim \frac{1}{T-t},$$ and which display some decoupling between the regular and the singular part of the solution at blow up time. We prove that this dynamic is unstable. More precisely, we show that any such solution with small super critical $L^2$ mass lies on the boundary of both $H^1$ open sets of global solutions that scatter forward and backwards in time, and solutions that blow up in finite time on the right in the log-log regime exhibited in \cite{MR1}, \cite{MR4}, \cite{R1}. We moreover exhibit some continuation properties of the scattering solution after blow up time and recover the chaotic phase behavior first exhibited in \cite{Mcpam} in the critical mass case.
\end{abstract}

\maketitle

%%%%%%%%%%%%%%%%%%%%%%%%%%%%
%%%%%%%%%%%%%%%%%%%%%%%%%%%%

\section{Introduction}

%%%%%%%%%%%%%%%%%%%%%%%%%%%%
%%%%%%%%%%%%%%%%%%%%%%%%%%%%

%%%%%%%%%%%%%%%%%%%%%%%%%%%%
%%%%%%%%%%%%%%%%%%%%%%%%%%%%

\subsection{Setting of the problem}

%%%%%%%%%%%%%%%%%%%%%%%%%%%%
%%%%%%%%%%%%%%%%%%%%%%%%%%%%

We consider in this paper the cubic two dimensional focusing nonlinear Schr\"odinger equation:
\be
\label{nls}
(NLS) \ \  \left   \{ \begin{array}{ll}
         iu_t=-\Delta u-|u|^{2}u, \ \ (t,x)\in \R\times \R^2,\\
         u(t_0,x)=u_0(x), \ \ u_0:\R^2\to \C, \ \ t_0\in \Bbb R.
         \end{array}
\right.
\ee
This is the special case of physical relevance of the $N$ dimensional $L^2$ critical (NLS):
\be
\label{nlsn}
  iu_t=-\Delta u-|u|^{\frac 4 N}u, \ \ (t,x)\in \R\times \R^N.
  \ee
From a result of Ginibre and Velo \cite{GV}, (\ref{nls}) is locally well-posed in $H^1=H^1(\R^2)$ and thus, for $u_0\in H^1$, there exists $t_0<T\leq +\infty$ and a unique solution $u(t)\in {\mathcal{C}}([t_0,T),H^1)$ to (\ref{nls}) and either $T=+\infty$, we say the solution is global, or $T<+\infty$ and then $\lim_{t\uparrow T}\|\nabla u(t)\|_{L^2}=+\infty$, we say the solution blows up in finite time. The Cauchy problem can also be solved in the critical $L^2$ space, \cite{CW}, in which case from standard Strichartz bound \cite{Str}, \cite{Cazenave}, finite time blow up is equivalent to $$\|u\|_{L^4([0,T),L^4_x)}=+\infty.$$ Note also that the equation is time reversible and the flow can similarly be solved backwards in time.\\
(\ref{nls}) admits the following conservation laws in the energy space $H^1$: 
$$
\left   . \begin{array}{lll}
         L^2-\mbox{norm}: \ \ \|u(t)\|^2_{L^2}=\|u_0\|_{L^2}^2;\\
         \mbox{Energy}:\ \ E(u(t,x))=\frac{1}{2}\int|\nabla u(t,x)|^2dx-\frac{1}{4}\int |u(t,x)|^{4}dx=E(u_0);\\
         \mbox{Momentum}:\ \ Im(\int\nabla u(t,x)\overline{u(t,x)}dx)=Im(\int\nabla u_0(x)\overline{u_0(x)}dx).
         \end{array}
\right .
$$
A large group of $H^1$ symmetries leaves the equation invariant: if $u(t,x)$ solves \fref{nls}, then $\forall (\lambda_0,\tau_0,x_0,\beta_0,\gamma_0)\in \RR_*^+\times\RR\times\RR^2\times\RR^2\times \RR$, so does 
\be
\label{symmetrygroup}
v(t,x)=\lambda_0u(\lambda_0^2t+\tau_0,\lambda_0 x+x_0-\beta_0t)e^{i\frac{\beta_0}{2}\cdot(x-\frac{\beta_0}{2} t)}e^{i\gamma_0}.
\ee
The scaling symmetry $v(t,x)=\lambda_0u(\lambda^2_0t,\lambda_0 x)$ leaves the $L^2$ space invariant and hence the problem is $L^2$ critical. The additional pseudo conformal symmetry 
\be
\label{vbvbeojr}
v(\tau,x)=\mathcal Cu(\tau,x)=\frac{1}{|\tau|}u\left(-\frac{1}{\tau},\frac{x}{|\tau|}\right)e^{i\frac{|x|^2}{4\tau}},
\ee
is a continuous transformation for $\tau\neq 0$ in the virial space $$\Sigma=\{xu\in L^2\}\cap H^1.$$

Following \cite{GNN}, \cite{KW}, let now $Q$ be the unique $H^1$ nonzero positive radial solution to
\be\label{eqQ}
\Delta Q-Q+Q^3=0,
\ee
then the variational characterization of $Q$ ensures that initial data $u_0\in H^1$ with $\|u_0\|_{L^2}<\|Q\|_{L^2}$ yield global and bounded solutions, \cite{W1}.  Moreover, for $u_0\in \Sigma$ or $u_0\in H^1$ and radial, the solution scatters forward and backward in time, or equivalently from standard Strichartz estimates:
\be
\label{scatter}
\|u\|_{L^4(\Bbb R,L^4_x)}<+\infty,
\ee
see \cite{Tao1} and \cite{Do1}, \cite{Do2} for the defocusing case and references therein. At the critical mass threshold $\|u_0\|_{L^2}=\|Q\|_{L^2}$, two new dynamics occur: the solitary wave non dispersive dynamics $$u(t,x)=Q(x)e^{it},$$ and the critical mass blow up solution generated by the pseudo conformal symmetry 
\be\label{nls:1}
S(t,x)=\frac{1}{|t|}\left(Qe^{it\frac{|y|^2}{4}}\right)\left(\frac{x}{|t|}\right)e^{-\frac{i}{t}}, \ \ t\in \Bbb R^*
\ee
 which blows up with the pseudo conformal speed $$\|\nabla u(t)\|_{L^2}\sim \frac{1}{|t|}\ \ \mbox{as} \ \ t\to 0.$$ From \cite{M1}, $S(t)$ is the unique up to the symmetries critical mass blow up solution. Through the pseudo conformal symmetry, this also classifies the solitary wave as the unique non dispersive critical mass solution in $\Sigma$ up to the symmetries, see \cite{visanetal} \cite{LiZh} for an extension to lower regularity.
%%%%%%%%%%%%%%%%%%%%%%%%%%%%
%%%%%%%%%%%%%%%%%%%%%%%%%%%%

\subsection{Three known dynamics}

%%%%%%%%%%%%%%%%%%%%%%%%%%%%
%%%%%%%%%%%%%%%%%%%%%%%%%%%%

We focus from now on and for the rest of this paper onto $H^1$ solutions with mass slightly above the critical one: 
\be
\label{smallnessltwo}
\|Q\|_{L^2}<\|u\|_{L^2}<\|Q\|_{L^2}+\alpha^*
\ee for some small enough universal constant $\alpha^*>0$. A general open problem is to classify the possible dynamics of the flow near the solitary wave $Q$. Three kind of regimes have been exhibited so far:\\

{\it (i) Scattering solutions}: Solutions that scatter forward and backward in time in the sense of \fref{scatter} have been exhibited with arbitrary mass, see for example \cite{Planchon}, and from standard Strichartz estimates \cite{Str}, see also \cite{Cazenave}, the corresponding set of initial data is open in $H^1$.\\

{\it (ii) Soliton like solutions and pseudo conformal blow up}: Solutions that scatter forward to $Q$ with super critical mass $$u(t)-Q-e^{it\Delta}u_{\infty}\to 0\ \ \mbox{in} \ \ H^1 \ \ \mbox{as} \ \ t\to+\infty$$ are constructed in the pioneering work by Bourgain and Wang \cite{BW} in dimensions $N=1, 2$, and in the further extension by Krieger and Schlag \cite{KS1} in dimension $N=1$. Through the pseudo conformal transformation \fref{vbvbeojr}, they equivalently correspond to super critical mass finite time blow up solutions with exact $S(t)$ blow up profile. 

\begin{Thm}[Bourgain, Wang \cite{BW}]
\label{thmbw}
Let $A_0$ be a given integer. Let $A\geq A_0$  a large enough integer. Let $$z^*\in X_A=\{f\in H^A \ \ \mbox{with} \ \ (1+|x|^A)f\in L^2\},$$ and let $z\in \mathcal C((T^*,0], X^A)$ be the solution to 
\be
\label{nlszstar}
\left\{\begin{array}{ll} i\pa_t z+\Delta z+z |z|^2=0,\\
z_{|t=0}=z^*,\end {array}\right.
\ee
where $T^*<0$ is the maximal time of existence of $z$. Assume that $z^*$ vanishes to high order at the origin:  
\be
\label{vanishinghigh}
D^{\alpha}z^*(0)=0 \ \ \mbox{for} \ \ |\alpha|\leq A-1.
\ee
Then for all $\theta\in \Bbb R$, there exists $t_0<0$ and a unique solution to $u^{\theta}_{BW}\in \mathcal C([t_0,0),X_{A_0})$ to \fref{nls} with 
\be
\label{esturestebourgianwang}
\|u^{\theta}_{BW}(t)-S(t)e^{i\theta}-z(t)\|_{X_{A_0}}\leq |t|^{A_0}.
\ee
\end{Thm}

For $t<0$, we will denote $$u_{BW}^{\theta}(t)$$ the Bourgain Wang solution which blows up at time $T=0$ with regular part $z^*$ and singular part $S(t)e^{i\theta}$ given by Theorem \ref{thmbw}. Note that $S(t)$ scatters to the left as $t\to -\infty$, and a further simple argument\footnote{Indeed, on can show that $u_{BW}^{\theta}$ satisfies \eqref{mardi}, and then Lemma \ref{lemmaopen} yields the result.} ensures that under the additional smallness assumption
\be
\label{smallness}
\|z^*\|_{H^A}<\alpha^*\ll 1,
\ee 
$u_{BW}^{\theta}$ scatters to the left in time:
\be
\label{scatterongbw}
u_{BW}^{\theta}\in \mathcal C((-\infty,0),H^{A_0}), \ \ \|u_{BW}^{\theta}(t,x)\|_{L^4((-\infty,-1],L^4)}<\infty.
\ee
Moreover, it blows up at $t=0$ with blow up speed given by the pseudo-conformal speed: 
\be
\label{speedconformal}
\|\nabla u_{BW}^{\theta}(t)\|_{L^2}\sim \frac{1}{|t|}.
\ee
From the time reversibility of the equation, we let for $t>0$ $u_{BW}^{\theta}(t)$ be the solution which blows up backwards at $t=0$ with $S(t)e^{i\theta}$ singular part and $z^*$ regular part, and $u_{BW}^{\theta}$ is global and scatters at $+\infty$ and blows up at the origin with speed \fref{speedconformal}.\\
In general, solutions blowing up with the pseudo conformal blow up speed \fref{speedconformal} have been conjectured to be {\it unstable} and to live on a codimension one manifold, see Krieger, Schlag \cite{KS1} for further results in this direction.\\

{\it (iii) Log-log blow up solutions}: Eventually, after the pioneering work by Perelman \cite{P} in dimension $N=1$, the existence of an $H^1$ open set of initial data leading to finite time blow up solutions in the so called log-log regime is proved in the series of works by Merle and Rapha\"el \cite{MR1}, \cite{MR2}, \cite{R1}, \cite{MR3}, \cite{MR4}, \cite{MR5}. These solutions concentrate a universal bubble of mass at blow up $$u(t,x)-\frac{1}{\lambda(t)}Q\left(\frac{x-x(t)}{\lambda(t)}\right)e^{i\gamma(t)}\to u^* \ \ \mbox{in} \ \ L^2 \ \ \mbox{as}\ \ t\to T=T(u)<+\infty $$ for some parameters $x(t)\to x(T)\in \Bbb R^2$, $\gamma(t)\to +\infty$ and a blow up speed 
\be
\label{loglogregime}
\|\nabla u(t)\|_{L^2}= \frac{\|\nabla Q\|_{L^2}}{\lambda(t)}, \ \ \lambda(t)=\sqrt{\frac{2\pi(T-t)}{\log |\log(T-t)|}}(1+o(1)) \ \ \mbox{as} \ \ t\to T.
\ee
\medskip 
Note that if we normalize the blow up time at $T=0$ and consider the pseudo conformal transformation \fref{vbvbeojr}, then $\mathcal Cu(\tau,x)$ is global on the right in time and behaves like a spreading bubble with anomalously slow decay: 
\be
\label{anomlousscattering}
\|\mathcal Cu(\tau,x)\|^4_{L^4}\sim \frac{\log|\log\tau|}{\tau} \ \ \mbox{as} \ \ \tau\to+\infty.
\ee
Note that this dynamic is by construction unstable by loglog blow up.\\

Finally, let us observe that the critical mass blow up solution $S(t)$ lies on the boundary of both $H^1$ open sets of global solutions that scatter forward and backwards in time and solutions that blow up in finite time on the right in the log-log regime due to the following explicit deformations:
\begin{itemize}
\item $u_{\eta}(-1)=(1-\eta)S(-1)$, $\eta>0$, has subcritical mass and thus the corresponding solution is global and scatters on both sides in time,  \cite{Tao1};
\item let $u_{\eta}(-1)=(1-\eta)S(-1)$, $\eta<0$ small,  and $v_{\eta}=\mathcal Cu$ its pseudo conformal transformation given by \fref{vbvbeojr}, then $v_{\eta}(1)=(1-\eta)Q$ has small super critical mass and $E(v_{\eta})<0$ from direct check, hence it blows up in finite time $1<\tau_\eta<+\infty$ in the log log regime, \cite{MR1}, \cite{MR2}, \cite{MR3}, and thus from \fref{vbvbeojr}, $u_\eta$ also blows up in finite time $T_{\eta}<0$ in the log log regime. 
\end{itemize}

%%%%%%%%%%%%%%%%%%%%%%%%%%%%
%%%%%%%%%%%%%%%%%%%%%%%%%%%%

\subsection{Statement of the result}

%%%%%%%%%%%%%%%%%%%%%%%%%%%%
%%%%%%%%%%%%%%%%%%%%%%%%%%%%

Our aim in this paper is to make some progress towards the understanding of the flow near $Q$ and the proof of the conjectured instability of the pseudo conformal type of blow up \fref{speedconformal}. We prove this instability for the special class of Bourgain Wang solutions. The particularity of these solutions is the decoupling in space between the singular $S(t)$ blow up bubble and the residual smooth $z^*$ part induced by the high degeneracy at blow up point \fref{vanishinghigh}. More precisely, we claim that Bourgain Wang solutions like the critical mass solution $S(t)$ lie on the border of both the $H^1$ open set of forward and backward scattering solutions, and the $H^1$ open set of solutions which blow up in finite time in the log-log regime \fref{loglogregime}. The following theorem is the main result of this paper:

\begin{Thm}[Strong instability of Bourgain-Wang solutions]
\label{thmmain}
Let $\a^*>0$ be a small enough universal constant and $A$ be a large enough integer. Let $z^*\in X_A$ radially symmetric satisfying the smallness assumption \fref{smallness} and the degeneracy at blow up point \fref{vanishinghigh}. Let $u_{BW}^0\in \mathcal C((-\infty,0),\Sigma)$ be the corresponding Bourgain-Wang solution given by Theorem \ref{thmbw} with $\theta=0$. Then there exists a continuous map $$\Gamma:[-1,1]\to \Sigma $$ such that the following holds true. Given $\eta\in [-1,1]$, let $u_{\eta}(t)$ be the solution to \fref{nls} with data $u_{\eta}(-1)=\Gamma(\eta)$, then:
\begin{itemize}
\item $\Gamma(0)=u_{BW}^0(-1)$ ie $\forall t<0$, $u_{\eta=0}(t)=u_{BW}^0(t)$ is the Bourgain Wang solution on $(-\infty,0)$ with blow up profile $S(t)$ and regular part $z^*$;
\item $\forall \eta\in(0,1]$, $u_{\eta}\in \mathcal C(\Bbb R,\Sigma)$ is global in time and scatters forward and backwards in the sense of \fref{scatter};
\item $\forall \eta\in [-1,0)$, $u_{\eta}\in \mathcal C((-\infty,T^*_{\eta}),\Sigma)$ scatters to the left and blows up in finite time  $T^*_{\eta}<0$ on the right {\it in the log-log regime} \fref{loglogregime} with
\be
\label{fitneitiem}
T^*_\eta\to 0\ \ \mbox{as}\ \ \eta\to 0.
\ee 
\end{itemize}
\end{Thm}

Theorem \ref{thmmain} will follow trough a complete description of the curve $\Gamma$ in terms of {\it explicit} modified profiles $\tilde{Q}_{\eta}(t)$ which are deformations of $S(t)$ {\it across} the pseudo conformal regime, see \fref{solutionexplicite} for precise definitions. In the log log regime $\eta<0$, the dynamics splits into three parts: the scattering dynamics for $t<-1$, the pseudo conformal blow up where the solution remains near $u_{BW}$ for $t<(1+\delta)T_{\eta}^*$, $0<\delta\ll1$, and eventually the log-log blow up for $t\sim T_{\eta}^*$. The control of the solution in each regime requires a specific analysis: Strichartz control at critical $L^2$ regularity for $t<-1$, a first monotonicity formula valid in the pseudo conformal regime\footnote{see Lemma \ref{lemma:timederivative}} only, and then the log log machinery based on {\it another} monotonicity formula from \cite{MR1}, \cite{MR2}, \cite{MR3}, \cite{MR4}, \cite{R1}.\\

In the global scattering regime $\eta>0$, an important outcome of the proof is that we may follow the flow of $u_{\eta}(t)$ for $\eta>0$ past the blow up time $t=0$ of the limiting Bourgain Wang solution as $\eta\to 0$. Indeed, a straightforward continuity argument ensures a strong limit before blow up time: $$\forall t<0,\ \ 
u_{\eta}(t)\to u^0_{BW}(t)\ \ \mbox{in}\ \ \Sigma\ \  \mbox{as}\ \ \eta\to 0.$$ The important problem is to understand what happens for $t>0$, and we recover the chaotic phase behavior first exhibited by Merle \cite{Mcpam} in the critical mass case:

\begin{Thm}[Continuation of $u_{\eta}$ after blow up time]
\label{themcontinuation}
Let $0<\eta\leq 1$ and $u_{\eta}\in \mathcal C(\Bbb R,\Sigma)$ be the global scattering solution built in Theorem \ref{thmmain}. Then:\\
(i) {\em Identification of the limit points}: Let $t>0$, then the limit points of the sequence $(u_{\eta}(t))_{\eta>0}$ as $\eta\to0$ are given by $$u_{BW}^{\theta}(t), \ \ \theta\in \Bbb R.$$
(ii) {\em Existence of converging subsequences}: $\forall \theta\in \Bbb R$, there exists a sequence $\eta_n\to 0$ such that 
 \be
 \label{opjepoueu}
 \forall t>0,\ \ u_{\eta_n}(t)\to u_{BW}^\theta(t) \ \ \mbox{in}\ \ \Sigma\ \  \mbox{as} \ \ n\to +\infty.
 \ee
\end{Thm}

In other words, the regularization after blow up time of the pseudo conformal Bourgain Wang blow up corresponds to a defocalization onto another Bourgain Wang profile with shifted blow up profile $S(t)e^{i\theta}$ and same regular part $z^*$, and where the phase shift displays a chaotic behavior as $\eta\to 0$ for $t>0$.\\

Eventually, the pseudo conformal symmetry \fref{vbvbeojr} yields the following corollary of instability of the manifold of forward scattering to $Q$.

\begin{Cor}[Instability of the manifold of forward scattering to $Q$]
\label{corforward}
Let $z^*, u_{\eta}$ as in the hypothesis of Theorem \ref{thmmain}. Then there exist continuous maps $$\Gamma^i:[-1,1]\to \Sigma, \ \ i=1,2 $$ such that the solution to \fref{nls} with data $v^{i}_{\eta}(1)=\Gamma^i(\eta)$ satisfies the following:
\begin{itemize}
\item $v^i_{\eta=0}\in \mathcal C(\Bbb R^+_*,\Sigma)$ and scatters to $Q$ to the right: $$v^i_0(\tau)-Q-\mathcal Cz(\tau)\to 0 \  \mbox{in} \ \ L^2\ \ \mbox{as} \ \ \tau\to +\infty.$$
\item $\forall \eta\in(0,1]$, $v^i_{\eta}\in \mathcal C(\Bbb R^*,\Sigma)$ scatters forward as $\tau\to+\infty$ in the sense of \fref{scatter}.
\item $\forall \eta\in [-1,0)$, $v^1_{\eta}\in \mathcal C((0,T^*_{\eta}),\Sigma)$  blows up in finite time  $0<T^*_{\eta}<+\infty$ on the right {\it in the log-log regime} with $T_{\eta}^*\to +\infty$ as $\eta\to 0$.
\item $\forall \eta\in [-1,0)$, $v^2_{\eta}\in \mathcal C((0,+\infty),\Sigma)$  and spreads to the right with the anomalous decay rate \fref{anomlousscattering}.
\end{itemize}
\end{Cor}

{\it Comments on the result}\\

{\it 1. Comparison with previous works:} The instability result of Theorem \ref{thmmain} and Corollary \ref{corforward} should be compared with the interesting classification results recently obtained by Nakanishi and Schlag \cite{NS1}, \cite{NS2}. These works lies in the continuation of former and more recent works which address the question of the existence and uniqueness of critical dynamics \cite{Mcpam}, \cite{M1}, \cite{DM1}, \cite{DM2}, \cite{DR}, \cite{RS} and the understanding of the center manifold near the ground state in super critical regimes \cite{S1}, \cite{KS1}, \cite{KS2}, \cite{KS3}, \cite{Be}. In \cite{NS2}, the dynamics of the flow near the ground state solitary wave for the $L^2$ super critical cubic (NLS) in dimension 3 are completely classified through a series of 9 different scenario which in particular identify the stable manifold of forward scattering to $Q$ as a threshold dynamic between the stable dynamics of finite time blow up and scattering. Note in passing that the scenario of anomalous scattering \fref{anomlousscattering} does not occur in \cite{NS2}. In the work \cite{NS2} however, finite time blow up is obtained through a virial argument without any qualitative information on the attained blow up regime. Moreover, the analysis by Krieger and Schlag \cite{KS1} suggests that despite some recent important progress on the construction and the understanding of the scattering manifold to $Q$ in super critical settings, it seems substantially more delicate to implement this kind of analysis to the $L^2$ critical case under consideration which is degenerate and critical with respect to a conservation law.\\

{\it 2. On the existence and instability of threshold dynamics}: The construction of the deformation $\Gamma$ in Theorem \ref{thmmain} will not rely on a fixed point like argument  as in \cite{BW}, \cite{S1}, \cite{KS2}, \cite{KS3}, \cite{KST1}, \cite{KST2}, \cite{KST3} but rather on softer compactness arguments as in \cite{Mc}, \cite{martel}, \cite{KS1}, \cite{RS}, \cite{HR}. Moreover, the threshold Bourgain Wang regime will be achieved as the limit from above and below of a family of solutions which each correspond to a {\it stable} dynamics, respectively scattering and log log blow up, and does not need to be constructed a priori. A similar strategy which gives together the existence and the instability of the threshold dynamics was used by Matano and Merle in their study of the super critical heat equation, see \cite{MM} and references therein. The extension of this strategy to the full family of pseudo conformal blow up solutions build in \cite{KS1} is open.\\

{\it 3. Extensions}: We have decided to work for simplicity in dimension $N=2$ and for radial data, but similar issues could be addressed for the $N$ dimensional $L^2$ critical (NLS) \fref{nlsn}. The extension to dimension $N=1$ and for non radial Bourgain Wang solution in dimension $N=1,2$ is straightforward\footnote{The non radial case would follow by introducing modulation on the translation parameter and Galilean deformation of the blow up profile which are lower order deformations, see \cite{MR6}, \cite{PR}, \cite{RS} for the full treatment of these issues in a similar setting.}. For higher dimensions, the roughness of the nonlinearity becomes a serious trouble. We will however along the proof of Theorem revisit the Bourgain Wang machinery and our new strategy has two main outcomes. The first one is that we will completely avoid the use of the pseudo conformal symmetry as an explicit symmetry and use instead a number of tools developed by Rapha\"el and Szeftel \cite{RS} for the study of minimal mass conformal blow up for an inhomogeneous (NLS) equation which yield a direct dynamical investigation of the pseudo conformal blow up without any need for a symmetry. In this sense, we expect the results of Theorem \ref{thmmain} and Theorem \ref{themcontinuation} to extend to a large class of nonlinear critical dispersive problems. The second outcome is that we will seriously diminish the regularity needed for the Bourgain Wang analysis, see in particular Remark \ref{remarkregularite}, which makes a partial extension of the above results in particular to dimension $N=3$ more reasonable to attain. This would require a separate analysis which remains to be done. 
Eventually, the smallness assumption \fref{smallness} on $z^*$ is used only to ensure global scattering behavior away from blow up. Without smallness,  similar results can be proved on the time interval $(-T^*_1,T^*_2)$ of existence of the solution to \fref{nlszstar} with very similar proofs and using in particular the large $L^2$ mass loglog analysis  in \cite{MR6}.

%%%%%%%%%%%%%%%%%%%%%%%%%%%%%%%%% 
%%%%%%%%%%%%%%%%%%%%%%%%%%%%%%%%%

\subsection{Strategy of the proof}

%%%%%%%%%%%%%%%%%%%%%%%%%%%%%%%%%
%%%%%%%%%%%%%%%%%%%%%%%%%%%%%%%%%

Let us give a brief insight into the strategy of the proof of Theorem \ref{thmmain} which follows from two main steps.\\

{\bf step 1} Approximate profiles.\\

The first step is to construct approximate profiles $\tilde{Q}_{\eta}(t)$ which are transversal to the pseudo conformal manifold. The instability direction is remarkably enough completely explicit and corresponds to profiles: 
$$\qte(t,x)=\frac{1}{\lambda_{\eta}(t)}\left[P_{\eta}(y)e^{-ib_{\eta}(t)\frac{|y|^2}{4}}\right]\left(\frac{x}{\lambda_{\eta}(t)}\right)e^{i\gamma_{\eta}(t)},
$$
for the explicit choice of modulation parameters:
\be
\label{parametersinitiallyonebis}
\lambda_{\eta}(t)=\sqrt{\eta+t^2}, \ \ b_{\eta}(t)=-t, \ \ \gamma_{\eta}(t)=1+\int_{-1}^t\frac{d\tau}{\eta+\tau^2},
\ee
and where $P_{\eta}$ is a suitable small deformation of $Q$. The profiles $P_{\eta}$ yield exact well localized scattering solutions for $\eta<0$, and $P_{\eta=0}=Q$ so that $$\tilde{Q}_{\eta=0}(t)=S(t).$$ For $\eta>0$, the self similar profile $P_{\eta}$ is unbounded in $L^2$ and the corresponding solution will blow up in finite time $T_{\eta}<0$ in the log-log regime instead of the self similar regime $\lambda(t)\sim\sqrt{\tilde{T}_{\eta}-t}$ with $\tilde{T}_{\eta}=-\sqrt{|\eta|}$ formally predicted by \fref{parametersinitiallyonebis}.\\

{\bf step 2} Backwards bootstrap of the flow.\\

We now pick an asymptotic smooth profile $z^*$ satisfying the degeneracy \fref{vanishinghigh}. We also assume the smallness condition \fref{smallnessltwo} which ensures that the solution $z$ to \fref{nlszstar} is global in time. Our analysis is based on a backwards integration of the flow in the continuation of the works \cite{Mc}, \cite{martel}, \cite{RS}.\\

\noindent{\bf (i)} For $\eta>0$, we let an initial data at time $t=0$ be $$u_{\eta}(t=0)=\qte(0)+z^*.$$ We then make a decomposition of the flow 
\be
\label{ceoieofie}
u_{\eta}=\qte+z(t)+\ut_{\eta}(t)
\ee
and claim a uniform backwards control of $\ut_\eta$ all the way down to some time $t=t_1$, where $t_1<0$ is independent of $\eta$: 
\be
\label{cbeieie}
 \|\ut_{\eta}(t)\|_{L^{\infty}([t_1,0],H^1)}\ll 1.
 \ee 
 The same argument also yields a similar control on $[0,|t_1|]$, and now from a standard continuity argument, the solution $u_{\eta}$ is  global and scatters at both infinities because so do $S(t)$ and $z(t)$. This yields both the scattering part of Theorem \ref{thmmain} and the continuation result of Theorem \ref{themcontinuation} where the chaotic phase is a direct consequence of the explicit formula for the phase \fref{parametersinitiallyonebis}.\\
 The key is therefore the backwards control \fref{cbeieie}. For this, we observe that the decomposition of \fref{ceoieofie} together with the strong decoupling \fref{vanishinghigh} induces an {\it almost critical mass} dynamics for $\tu_{\eta}$. We are thus able to use the Lyapounov machinery exhibited in \cite{RS} in the pseudoconformal regime to control the flow backwards in time.\\

\noindent{\bf (ii)} For $\eta>0$, an explicit computation shows that the solution to \fref{nls} with data $$u_{\eta}(T^1_{\eta})=\qte(T^1_{\eta})+z^* \ \ \mbox{at} \ \ T^1_{\eta}\sim -\sqrt{\eta}$$ blows up in the log log regime forward in time. Here we use the sharp open characterization of the log log set exhibited in \cite{R1} which allows one to see directly loglog blow up on the initial data. We now need to run backwards the control of this solution. The key observation is that at the time $T^1_{\eta}$, the solution is still very much in the pseudoconformal regime and the backwards reintegration of the flow is completely similar to the one of \fref{cbeieie}.\\

\noindent{\bf (iii)} We now pass to the limit as $\eta\to 0$ and obtain in both cases a Bourgain Wang type solution with regular part $z^*$ as $t\to 0$. It remains to show some uniqueness statement about this type of dynamics  to ensure the continuity of the map $$\Gamma:\eta\mapsto u_{\eta}(-1)$$ at $\eta=0$. This is a simple but slightly technical claim which proof relies again on the Lyapounov machinery exhibited in \cite{RS}, and this is the part which requires a high order degeneracy in \fref{vanishinghigh}.\\
 
This paper is organized as follows. In section \ref{sectionone}, we build the deformation $\qte$ of the $S(t)$ blow up profile. In section \ref{backwrds}, we build the machinery of backwards control of the pseudo conformal flow. The proof of the main theorems is then concluded in section \ref{seciofheio}.\\

\noindent{\bf Acknowledgments.} F.M. is supported by ANR Projet Blanc OndeNonLin. P.R and J.S are supported by ANR jeunes chercheurs SWAP. \\

\noindent{\bf Notations:} We introduce the differential operator $$\Lambda f=f+y\cdot\nabla f \ \ \mbox{($L^2$ scaling)}.$$  Let $L=(L_+,L_-)$ the matrix linearized operator close to the ground state with:
\be
\label{deflpluslmoins}
L_+=-\Delta +1-3Q^2, \ \ L_-=-\Delta+1-Q^2.
\ee
We recall that $L$ restricted to radial functions has a generalized nullspace characterized by the following algebraic identities:
\be
\label{structurekernelpair}
L_-Q=0,\ \ L_+(\Lambda Q)=-2Q, \ \ L_-(|y|^2Q)=-4\Lambda Q, \ \ L_+\rho=|y|^2Q, 
\ee
where $\rho$ is the unique radial $H^1$ solution to 
\be
\label{defrho}
L_+\rho=|y|^2Q.
\ee

%%%%%%%%%%%%%%%%%%%%%%%%%%%%%%%%% 
%%%%%%%%%%%%%%%%%%%%%%%%%%%%%%%%%

\section{Approximate solutions}
\label{sectionone}

%%%%%%%%%%%%%%%%%%%%%%%%%%%%%%%%%
%%%%%%%%%%%%%%%%%%%%%%%%%%%%%%%%%

In this section, we construct the leading order terms in the dynamics we will consider. The regular part will be given by $z(t)$ solution to \fref{nlszstar} for which we recall some well known regularity statements, and the singular dynamics will be approximated by an {\it explicit} deformation of $S(t)$.

%%%%%%%%%%%%%%%%%%%%%%%%%%%%%%%%%
%%%%%%%%%%%%%%%%%%%%%%%%%%%%%%%%%

\subsection{Regular part of the dynamic}

%%%%%%%%%%%%%%%%%%%%%%%%%%%%%%%%%
%%%%%%%%%%%%%%%%%%%%%%%%%%%%%%%%%

From now on and for the rest of the paper, we let 
\be
\label{conditionm}
m\geq 3
\ee 
be an integer and consider $z^*\in H^{2m+3}$ radially symmetric with  
\be
\label{smallnessbis}
\|z^*\|_{H^{2m+3}}<\alpha^*
\ee for some small enough universal constant $\alpha^*>0$. We moreover assume the flatness at the origin:
\be
\label{flatattheorigin}
\forall 0\leq k\leq 2m, \ \ \frac{d^kz^*}{dr^k}|_{r=0}=0.
\ee
Let $z(t)$ be the solution to
\be
\label{defutilde}
\left\{\begin{array}{ll} i\pa_t z+\Delta z+z |z|^2=0,\\
z_{|t=0}=z^*.\end {array}\right .
\ee
then standard Cauchy theory coupled with the smoothness of the nonlinearity ensures:
\begin{Lemma}[Regularity of $z$]
\label{regularityz} 
We have $z\in \mathcal C(\Bbb R,H^{2m+3})$. Moreover:\\
(i) $z$ scatters ie 
\be
\label{zscatter}
\|z\|_{L^4(\Bbb R,L^4_x)}\lesssim\a^*;
\ee
(ii) there holds the local in time propagation of regularity: 
\be
\label{cnofeof} 
\|z\|_{L^{\infty}([-1,1],H^{2m+3})}\lesssim \a^*;
\ee
(iii) there holds the propagation of degeneracy at the origin: $\forall t\in [-1,1]$,
\be
\label{pointwiseboudnut}
 |z(t,x)|\lesssim \alpha^*(|t|+|x|^2)^{m+1},\, |\pa_rz(t,x)|\lesssim \alpha^*(|t|+|x|^2)^{m+\frac12},\, |\pa^2_rz(t,x)|\lesssim \alpha^*(|t|+|x|^2)^m.
\ee
\end{Lemma}

{\bf Proof of Lemma \ref{regularityz}}: The proof is standard and we briefly recall the main facts. \eqref{zscatter} follows from the global in time small data $L^2$ critical Cauchy theory, \cite{CW}. Next, the regularity of the cubic nonlinearity and Strichartz estimates \cite{Str} with the 2 dimensional admissible pair $(4,4)$ imply:
\bee
\|z\|_{L^{\infty}([-1,1],H^{2m+3})}&\lesssim&  \|z^*\|_{H^{2m+3}}+\||z|^2z\|_{L^{\frac{4}{3}}([-1,1],W^{2m+3,\frac{4}{3}}_x)}\\
&\lesssim& \a^*+ \|z\|_{L^{\infty}([-1,1],H^{2m+3}_x)}\|z\|^2_{L^{\infty}([-1,1],L^8_x)}\\
&\lesssim& \a^*+\|z\|_{L^{\infty}([-1,1],H^{2m+3}_x)}\|z\|^2_{L^{\infty}([-1,1],H^1_x)}
\eee
which ensures \eqref{cnofeof} for $\alpha^*$ small enough. Finally, \fref{flatattheorigin}, \fref{defutilde}, \fref{cnofeof} and a Taylor expansion at $(t,x)=(0,0)$ on $z$ yield \fref{pointwiseboudnut}. This concludes the proof of Lemma \ref{regularityz}.\\

Let us also recall the openness of the set of scattering solutions which follows from standard Strichartz estimates again, \cite{CW}:

\begin{Lemma}[The set of scattering solutions is open]
\label{lemmaopen}
Let $u\in \mathcal C([0+\infty),H^1)$ be a global solution to \fref{nls} which scatters on the right: $$\|u\|_{L^4([0,+\infty),L^4_x)}<+\infty.$$ Then there exists $\alpha>0$ such that for all $v_0\in H^1$ with $\|u_0-v_0\|_{L^2}<\alpha$, the solution $v$ to \fref{nls} with initial data $v_0$ is global $v\in \mathcal C([0+\infty),H^1)$ and scatters to the right $$\|v\|_{L^4([0,+\infty),L^4_x)}<+\infty.$$
\end{Lemma}

%%%%%%%%%%%%%%%%%%%%%%%%%%%%%%%%% 
%%%%%%%%%%%%%%%%%%%%%%%%%%%%%%%%%

\subsection{Approximate solution for the singular part of the dynamic}

%%%%%%%%%%%%%%%%%%%%%%%%%%%%%%%%%
%%%%%%%%%%%%%%%%%%%%%%%%%%%%%%%%%

Let us now construct the full family of self similar profiles which generalizes the construction in \cite{PLSS}, \cite{MR2}.\\
Let $u$ be a solution to \fref{nls} and let us pass to the self similar variables 
\be\label{mercredi}
u(t,x)=\frac{1}{\lambda(t)}v\left(t,\frac{x}{\lambda(t)}\right)e^{i\gamma(t)}, \ \ \frac{ds}{dt}=\frac{1}{\l^2},
\ee
which leads  to 
\be\label{mercredi1}
i\pa_sv+\Delta v-v+ib\Lambda v+v|v|^2=0, \ \gamma_s=1, \ \ -\lsl=b.
\ee
Let a parameter $\eta\in\R$. We look for specific solutions of the form $$v(s,y)=P_\eta(y)e^{-i\frac{b|y|^2}{4}}$$ which thus have to satisfy:
\be
\label{explicitsoltuions}
\Delta P_\eta-P_\eta+\frac{b_s+b^2}{4}|y|^2P_\eta+P_\eta|P_\eta|^2=0.
\ee
We now anticipate the law $b_s+b^2=-\eta$ and build the corresponding profiles from standard perturbative and elliptic techniques. 

\begin{Lemma}[Existence of well localized profiles]
\label{scatteringprofiles}
There exists a universal constant $\eta^*>0$ and a smooth map $\eta \to P_{\eta}\in H_{rad}^{\infty}$  defined on $(-\eta^*,\eta^*)$ with the following properties: $P_{\eta}$ is radially symmetric and real positive,  and solves an equation of the form
\be
\label{equationpeta}
\Delta P_{\eta}-P_{\eta}-\frac{\eta}{4}|y|^2P_{\eta}+P_{\eta}^3=-\Psi_{\eta}
\ee
where $\Psi_{\eta}$ is identically 0 in the case $0\leq \eta<\eta^*$, and some well localized and exponentially small error in the case $-\eta^*<\eta<0$: for all polynomial $q(y)$, $\forall k\geq 0$, 
\be\label{mercredi2}
\left\|q(y)\frac{d^k}{dy^k}\Psi_{\eta}(y)\right\|_{L^{\infty}}\lesssim e^{-\frac{C_{k,q}}{\sqrt{\eta}}}, \ \ \mbox{Supp}\Psi_{\eta}\subset\left\{\frac{1}{\sqrt{\eta}}, \frac{4}{\sqrt{\eta}}\right\}.
\ee
Moreover, for all $-\eta^*<\eta<\eta^*$, $P_{\eta}$ is a smooth function of $(\eta,y)$ with 
\be
\label{derivee}
\frac{\partial P_\eta}{\partial \eta}|_{\eta=0}=-\frac14\rho,
\ee
where $\rho$ is the unique $H^1$ radially symmetric solution to \eqref{defrho}, and there holds the uniform decay estimates: 
\be
\label{uniformdeacy}
\forall k\geq 0, \ \ \forall -\eta^*<\eta<\eta^*, \ \left|\frac{d^k}{dy^k}\frac{\partial P_{\eta}}{\partial \eta}\right|\leq e^{-c_k|y|}, \ \  \left|\frac{d^k}{dy^k}(P_{\eta}-Q)\right|\leq \eta e^{-c_k|y|}. 
\ee
\end{Lemma}

\begin{remark} The sign $\eta\geq 0$ in \fref{equationpeta} ensures the uniform ellipticity of the operator $\Delta -1-\frac{\eta|y|^2}{4}$ independently of $\eta$ and hence the uniformity in $\eta$ of the exponential localization bounds \fref{uniformdeacy}. The profiles $P_{\eta}$ for $\eta<0$ correspond to the almost self similar solutions build in \cite{MR2}, \cite{MR4} to which we refer for further details.
\end{remark}

\begin{remark} The invariant of mass and energy can be computed explicitly in the limit $\eta\to 0$. Indeed, from \fref{defrho}, \eqref{derivee}, and standard computations: 
\be\label{invariantl2peta}
\frac{d}{d\eta}\left\{\|P_{\eta}\|_{L^2}^2\right\}|_{\eta=0}=-\frac{1}{2}(Q,\rho)=\frac{1}{4}(|y|^2Q,\Lambda Q)=-\frac{1}{4}\|yQ\|^2_{L^2}.
\ee 
For the energy, we multiply \fref{equationpeta} by the Pohozaev multiplier $\Lambda P_{\eta}$ and  conclude: 
\be
\label{coputationngry}
E(P_{\eta})=\frac{\eta}{8}\|yQ\|_{L^2}^2+O\left(\eta^2+e^{-\frac{c}{\sqrt{\eta}}}\right)=\frac{\eta}{8}\|yQ\|_{L^2}^2+o(\eta).
\ee
\end{remark}
Given small parameters $(b,\eta)$ we define:
\be
\label{defqetab}
Q_{\eta,b}=P_{\eta}e^{-i\frac{b|y|^2}{4}},
\ee 
which satisfies:
\be\label{eqqb}
\Delta Q_{\eta,b}-Q_{\eta,b}+ib\Lambda Q_{\eta,b}-(b^2+\eta)\frac{|y|^2}{4}Q_{\eta,b}+Q_{\eta,b}|Q_{\eta,b}|^2=-\Psi_{\eta,b},
\ee
where $\Psi_{\eta,b}=0$ if $\eta\geq 0$ and $\Psi_{\eta,b}=\Psi_\eta e^{-i\frac{b|y|^2}{4}}$ if $\eta<0$. Following \fref{mercredi} \eqref{mercredi1}, we now let
\be
\label{solutionexplicite}
\qte(t,x)=\frac{1}{\lambda_{\eta}(t)}Q_{\eta,b_\eta(t)}\left(\frac{x}{\lambda_{\eta}(t)}\right)e^{i\gamma_{\eta}(t)},
\ee
where the parameters solve the finite dimensional system: 
\be
\label{finitediemnsianlparamaters}
-\frac{(\lambda_{\eta})_s}{\lambda_{\eta}}=b_{\eta}, \ \ (\gamma_{\eta})_s=1, \ \ (b_{\eta})_s+b_{\eta}^2=-\eta
\ee 
for the explicit solution given by:
\be
\label{parametersinitiallyone}
\lambda_{\eta}(t)=\sqrt{\eta+t^2}, \ \ b_{\eta}(t)=-t, \ \ \gamma_{\eta}(t)=1+\int_{-1}^t\frac{d\tau}{\eta+\tau^2}.
\ee
Note that $\tilde{Q}_{\eta}$ transitions smoothly from an {\it exact} global scattering profile for $\eta>0$ to an {\it almost} self similar solution for $\eta< 0$ which would formally blow up at $t=-\sqrt{|\eta|}<0$, and through the $S(t)$ blow up solution which blows up in the conformal regime at $t=0$ and corresponds to $\eta=0$.
%%%%%%%%%%%%%%%%%%%%%%%%%%%%%%%%% 
%%%%%%%%%%%%%%%%%%%%%%%%%%%%%%%%%

 %%%%%%%%%%%%%%%%%%%%%%%%%%%%%%%%% 
%%%%%%%%%%%%%%%%%%%%%%%%%%%%%%%%%

\section{Backwards control of the pseudo conformal blow up}
\label{backwrds}

%%%%%%%%%%%%%%%%%%%%%%%%%%%%%%%%%
%%%%%%%%%%%%%%%%%%%%%%%%%%%%%%%%%

We develop in this section some analytical tools to control the (NLS) flow near $Q$ backwards in time and in the pseudo conformal regime. We will use the decoupling \fref{pointwiseboudnut} to reduce ourselves to an almost critical mass setting and use the mixed Energy/Morawetz Lyapounov type functional exhibited in \cite{RS}, and which yields an {\it increasing} Lyapounov function hence suitable to integrate the flow backwards from the singularity. As we explained in the introduction, the backwards integration of the flow after respectively the scattering time for $\eta>0$ or the log-log time for $\eta<0$ is very similar, and we deal with both cases simultaneously. The heart of our analysis is the following:

\begin{Prop}[Backwards control]
\label{propbackwardintegration}
Let 
\be
\label{conditionmprop}
m\geq 3
\ee 
and $z^*\in H^{2m+3}$ radially symmetric satisfying \fref{smallnessbis}, \fref{flatattheorigin} for some $\alpha^*>0$ small enough. Let $z(t)\in \mathcal C(\Bbb R,H^{2m+3})$ be the solution to \fref{defutilde}. Let $|\eta|<\eta^*(\alpha^*)$ small enough and a time $T_\eta\in[-1,0] $ such that:
\be\label{dataatteta0}
\sqrt{|\eta|}+|t|\lesssim \l_\eta(t)\textrm{ for all }t\in[-1,T_\eta].
\ee
Let an arbitrary sequence $$\eta\mapsto \gamma^0_\eta\in \Bbb R$$ and $u_\eta$ be the solution to \fref{nls} with data at $t=T_{\eta}$:
\be\label{dataatteta}
u_\eta(T_\eta)=\qte(T_\eta)e^{i\gamma^0_{\eta}}+z(T_\eta).
\ee
Then:\\
(i) {\em Local backwards control}: there exists a time $t_1<T_{\eta}$ independent of $\eta$ such that  $u_{\eta}\in \mathcal C([t_1,T_{\eta}],H^{\frac{3}{2}})$ and admits a decomposition of the form $$u_{\eta}(t,x)=\frac{1}{\lambda(t)}(Q_{\eta,b(t)}+\e)\left(t,\frac{x}{\lambda(t)}\right)e^{i(\gamma(t)+\gamma^0_{\eta})}+z(t,x)$$ with the uniform bounds in $\eta$: $\forall t\in [t_1,T_\eta]$,
\be
\label{uniformbounds}
\frac{1}{\lambda^{m-1}}\left|\frac{\lambda_\eta}{\lambda}-1\right|+\frac{|b_{\eta}-b|}{\lambda^m}+\frac{|\gamma-\gamma_{\eta}|}{\lambda^{m-2}}+\frac{\|\e\|_{H^1}}{\l^{m+1}}+\frac{\|\e\|_{H^{\frac 32}}}{\l^{m-\frac 32}}\lesssim \alpha^*\ee
where $b_\eta, \lambda_\eta, \gamma_\eta$ are given by \eqref{parametersinitiallyone}.\\
(ii) {\it Global backwards control}: $u_{\eta}\in \mathcal C((-\infty,T_{\eta}),H^{\frac 32})$ and scatters backwards: 
\be
\label{scaterrsbacjgfiards}
\|u_{\eta}\|_{L^4((-\infty,t_1),L^4_x)}\lesssim 1.
\ee
\end{Prop}

The proof of Proposition \ref{propbackwardintegration} relies on a bootstrap argument. Indeed, from standard modulation theory\footnote{see \cite{RS} for further details} and the smallness of $\lambda_{\eta}(T_{\eta})$, we may consider a time $t_2<T_\eta$ a priori depending on $\eta$ such that $\forall t\in [t_2,T_{\eta}]$:
\be
\label{decompopofupe}
u=\frac{1}{\lambda(t)}(Q_{\eta,b(t)}+\e_{z}+\e)\left(t,\frac{x}{\lambda(t)}\right)e^{i(\gamma(t)+\gamma_{\eta}^0)}=w+\tilde{u}
\ee 
with
\be
\label{defomega}
w(t,x)=\frac{1}{\lambda(t)}(Q_{\eta,b(t)}+\e_z)\left(t,\frac{x}{\lambda(t)}\right)e^{i(\gamma(t)+\gamma^0_{\eta})}, \ \ z(t,x)=\frac{1}{\lambda(t)}\e_z\left(t,\frac{x}{\lambda(t)}\right)e^{i(\gamma(t)+\gamma^0_{\eta})}.
\ee
The uniqueness of the decomposition \fref{decompopofupe} is ensured through the choice of orthogonality conditions:
\be
\label{orthee}
Re(\e,\overline{|y|^2Q_{\eta,b}})=0, \ \ Im(\e,\overline{\Lambda Q_{\eta,b}})=0, \ \  Im(\e,\overline{\tilde{\rho}})=0
\ee
where $$\tilde{\rho}=\rho e^{-i\frac{b|y|^2}{4}}, \ \ L_+\rho=|y|^2Q.$$
Moreover, $(\lambda,b,\gamma,\e)(T_{\eta})=(\lambda_{\eta},b_{\eta},\gamma_{\eta},0)(T_{\eta})$ from \fref{solutionexplicite}, \fref{dataatteta}. Thus, from a straightforward continuity argument, we may assume that $u$ satisfies the following a priori bounds: 
\be
\label{aprioribound}
\sup_{t\in[t_2,T_{\eta}]}\left\{\frac{\|\e\|_{H^1}}{\l^{m+1}}\right\}\leq K,\, \sup_{t\in[t_2,T_{\eta}]} \left\{\frac{1}{\lambda^{m-1}}\left|\frac{\lambda_\eta}{\lambda}-1\right|+\frac{|b_{\eta}-b|}{\lambda^m}\right\}\leq 1,
\ee
where $K$ is a small enough constant. We now claim the following bootstrap Lemma:

\begin{Lemma}[Backwards bootstrap bound on $\e$]
\label{bootlemma}
Assume that $K$ in \fref{aprioribound} has been chosen small enough -independent of $\eta$-. Then there exists a small time $t_1<0$ -independent of $\eta$- such that $u$ satisfies on $[t_1,T_{\eta}]$ the improved bound: 
\be
\label{estktun}
\sup_{t\in[t_1,T_{\eta}]} \left\{\frac{1}{\lambda^{m-1}}\left|\frac{\lambda_\eta}{\lambda}-1\right|+\frac{|b_{\eta}-b|}{\lambda^m}+\frac{|\gamma-\gamma_{\eta}|}{\lambda^{m-2}}+\frac{\|\e\|_{H^1}}{\l^{m+1}}+\frac{\|\e\|_{H^{\frac 32}}}{\l^{m-\frac 32}}\right\}\lesssim \sqrt{\alpha^*}.
\ee
\end{Lemma}

\begin{remark}
\label{remarkltwo} From the proof, the size of $t_1$ is uniform as $\alpha^*\to 0$. 
\end{remark}

Let us conclude the proof of Proposition \ref{propbackwardintegration} assuming Lemma \ref{bootlemma}.\\

{\bf Proof of  Proposition \ref{propbackwardintegration}}: The local backwards control \fref{uniformbounds} follows from \fref{estktun}. We now claim at $t=t_1$ the uniform control: $\forall |\eta|\leq \eta^*(\alpha^*)$ small enough,
\be
\label{inojeoefg}
\|e^{-i\gamma^0_\eta}u_{\eta}(t_1)-S(t_1)\|_{L^2}\lesssim \sqrt{\alpha^*}
\ee
which together with Remark \ref{remarkltwo}, Lemma \ref{lemmaopen} and the observation that $S(t)$ scatters as $t\to -\infty$ yields the global control \fref{scaterrsbacjgfiards}. This concludes the proof of Proposition \ref{propbackwardintegration} provided \eqref{inojeoefg} holds.

Now, we turn to the proof of \fref{inojeoefg}. Let 
$$\tq(t_1)=\frac{1}{\lambda(t_1)}Q_{\eta,b}\left(t_1,\frac{x}{\lambda(t_1)}\right)e^{i(\gamma(t_1)+\gamma_\eta^0)}.$$
Then, we have:
\bea
\label{domani}&&\|S(t_1)-e^{-i\gamma^0_\eta}\tq(t_1)\|_{L^2}\\
\nonumber&\lesssim& \|Q-P_\eta\|_{L^2}+|b+t|+\left|\frac{t}{\l}-1\right|+\left|\gamma-\frac{1}{t}\right|\\
\nonumber&\lesssim& \|Q-P_\eta\|_{L^2}+|b-b_\eta|+\left|\frac{\l_\eta}{\l}-1\right|+\left|\gamma-\gamma_\eta\right|+|b_\eta+t|+\left|\frac{t}{\l_\eta}-1\right|+\left|\gamma_\eta-\frac{1}{t}\right|\\
\nonumber&\lesssim& |\eta|+\sqrt{\a^*},
\eea
where we used \eqref{uniformdeacy}, \eqref{parametersinitiallyone} and \eqref{estktun} in the last inequality. Now, the decomposition \fref{decompopofupe} together with the bounds \eqref{smallnessbis}, \eqref{estktun} and \eqref{domani} yields:
\bee \|u_{\eta}(t_1)e^{-i\gamma^0_\eta}-S(t_1)\|_{L^2} & \lesssim & \|z(t_1)\|_{L^2}+\|\tilde{u}(t_1)\|_{L^2}+\|S(t_1)-e^{-i\gamma^0_\eta}\tq(t_1)\|_{L^2}\\
& \lesssim & |\eta|+\sqrt{\alpha^*},
\eee
which implies \fref{inojeoefg} for $|\eta|<\eta^*(\alpha^*)$ small enough. This concludes the proof of Proposition \ref{propbackwardintegration}.\\

The rest of this section is devoted to the proof of Lemma \ref{bootlemma} which relies on the Lyapounov functional approach developed in \cite{RS}.

%%%%%%%%%%%%%%%%%%%%%%%%%%%%%%%%% 
%%%%%%%%%%%%%%%%%%%%%%%%%%%%%%%%%

\subsection{Modulation equations}

%%%%%%%%%%%%%%%%%%%%%%%%%%%%%%%%%
%%%%%%%%%%%%%%%%%%%%%%%%%%%%%%%%%

Let us introduced the rescaled time $$\frac{ds}{dt}=\frac{1}{\lambda^2}.$$ Let $w$ be the refined profile given by \fref{defomega} and $W$ its renormalized version:  
\bea
\label{defqtilde}
w(t,x)=\tq+z=\frac{1}{\lambda(t)}W\left(t,\frac{x}{\lambda(t)}\right)e^{i(\gamma(t)+\gamma^0_{\eta})},\ \ W(s,y)=Q_{\eta,b(s)}(y)+\e_{z}(s,y),\\  
\nonumber\tq=\frac{1}{\lambda(t)}Q_{\eta,b}\left(t,\frac{x}{\lambda(t)}\right)e^{i(\gamma(t)+\gamma_\eta^0)}, \ \ z=\frac{1}{\lambda(t)}\e_{z}\left(t,\frac{x}{\lambda(t)}\right)e^{i(\gamma(t)+\gamma^0_{\eta})}.
\eea
From Lemma \ref{scatteringprofiles} and \fref{cnofeof}, we have the bound:
\be
\label{aprioriboundwsw}
\|w\|_{L^2}\lesssim 1, \ \ \|\nabla w\|_{L^2}\lesssim \frac{1}{\lambda}, \ \ \|w\|_{H^{\frac{3}{2}}}\lesssim \frac{1}{\lambda^{\frac{3}{2}}}.
\ee
Moreover, $w$ satisfies the equation
\be
\label{eqwgobale}
i\pa_tw+\Delta w+w|w|^2=\psi=\frac{1}{\lambda^3}\Psi\left(t,\frac{x}{\lambda(t)}\right)e^{i(\gamma(t)+\gamma^0_\eta)},
\ee
\be
\label{defpsiacaculaer}
\Psi =  -i\left(\lsl+b\right)\Lambda Q_{\eta,b}+\frac14(b_s+b^2+\eta)|y|^2Q_{\eta,b}-\tgamma_sQ_{b,\eta}-\Psi_{\eta,b}+\tilde{R},
\ee
\be
\label{defrtide}
\tilde{R}=(Q_{\eta,b}+\e_z)|Q_{b,\eta}+\e_z|^2-Q_{\eta,b}|Q_{\eta,b}|^2-\e_z|\e_z|^2,
\ee
and where $\tgamma_s=\gamma_s-1.$ We then decompose $u=w+\tilde{u}$ so that $\tilde{u}$ satisfies: 
\be\label{un3}
i\p_t\tu+\Delta\tu+(|u|^2u-|w|^2w)=-\psi.
\ee
We rewrite the interaction term as:
$$
|u|^2u-|w|^2w =  \frac{1}{\lambda^3}\big[2\e|Q_{\eta,b}|^2+Q_{\eta,b}^2\overline{\e}+R(\e)\big]\left(t,\frac{x}{\lambda(t)}\right)e^{i(\gamma(t)+\gamma^0_\eta)}\\
$$
with 
\be
\label{defrepsilon}
R(\e)=2\left(|W|^2-|Q_{\eta,b}|^2\right)\e+\left(W^2-Q_{\eta,b}^2\right)\overline{\e}+  \left(2W|\e|^2+\overline{W}\e^2+\e|\e|^2\right).
\ee
We then decompose $$\ \ Q_{\eta,b}=\S+i\T, \ \ \e=\e_1+i\e_2, \ \ R(\e)=R_1(\e)+iR_2(\e), \ \ \tilde{R}=\tilde{R}_1+i\tilde{R}_2$$ in terms of real and imaginary part and obtain the equation satisfied by $\e$:
\bea
\label{eqeonereal}
\nonumber \pa_s\e_1-M_2(\e)+b\Lambda \e_1 & = & \left(\lsl+b\right)(\Lambda \S+\Lambda\e_1)-\frac14(b_s+b^2+\eta)|y|^2\T+\tgamma_s(\T+\e_2)\\
& + & Im(\Psi_{\eta,b})-\tilde{R}_2-R_2(\e),
\eea
\bea
\label{eqeoneim}
\nonumber
\pa_s\e_2+M_1(\e)+b\Lambda \e_2 & = & \left(\lsl+b\right)(\Lambda \T+\Lambda\e_2)+\frac14(b_s+b^2+\eta)|y|^2\S-\tgamma_s(\S+\e_1)\\
& - & Re(\Psi_{\eta,b})+\tilde{R}_1+R_1(\e),
\eea
where $(M_1,M_2)$ are small deformations of the linearized operators $(L_+,L_-)$ close to $Q$:
$$M_1(\e)=-\Delta \e_1+\e_1-(3\S^2+\T^2)\e_1-2\S\T\e_2,  \ \ M_2(\e)=-\Delta \e_2+\e_2-(3\T^2+\S^2)\e_2+2\S\T\e_1.$$ 
We now claim the following preliminary estimates on the decomposition which are a consequence of the orthogonality conditions \fref{orthee}:

\begin{Lemma}[Preliminary estimates on the decomposition]
\label{lemmapreliminary}
 (i) Degeneracy of the unstable direction: there holds
\be
\label{degeneracyeunqlemma}
|(\e_1,Q)|\lesssim \alpha^*\lambda^{m+1}+K\lambda^{m+2}.
\ee
(ii) Modulation equations: Let $$Mod(t)=\left(b_s+b^2+\eta,\lsl+b,\tgamma_s\right)$$ then 
\be
\label{controlmodulation}
|Mod(t)|\lesssim \alpha^*\lambda^{m+1}+K\lambda^{m+2},
\ee
and 
\be
\label{controlmodulationbis}
\left|\lsl+b\right|\lesssim \alpha^*\lambda^{m+2}+K\lambda^{m+2}.
\ee
\end{Lemma}

{\bf Proof of Lemma \ref{lemmapreliminary}}\\

We compute the modulation equations by taking the inner product of \fref{eqeonereal}, \fref{eqeoneim} with the well localized in space directions corresponding to the orthogonality conditions \fref{orthee}.\\

{\bf step 1} Inner products.\\

We compute the inner products needed to obtain \eqref{degeneracyeunqlemma}, \eqref{controlmodulation}, and \eqref{controlmodulationbis}. This computation has been made in \cite{RS} near $P_{\eta=0}=Q$ and the same computation is valid up to $O(\eta\|\e\|_{L^2})$ terms, and leads to:
\be\label{pscal0}
 \left(-M_2(\e)+b\Lambda\e_1,\S\right)+\left(M_1(\e)+b\Lambda\e_2,\T\right)=  O(\l^2\|\e\|_{L^2}),
\ee
\be\label{pscal1}
-\left(-M_2(\e)+b\L\e_1,\Lambda\T\right)+\left(M_1(\e)+b\Lambda\e_2,\Lambda\S\right)=  -2\Re(\e,\overline{Q_{\eta,b}})+O(\l^2\|\e\|_{L^2}),
\ee
\be\label{pscal2}
 \left(-M_2(\e)+b\Lambda\e_1,|y|^2\S\right)+\left(M_1(\e)+b\Lambda\e_2,|y|^2\T\right)=  O(\l^2\|\e\|_{L^2}),
\ee
\be\label{pscal3}
 -\left(-M_2(\e)+b\Lambda\e_1,\rho_2\right)+\left(M_1(\e)+b\Lambda\e_2,\rho_1\right)=  O(\l^2\|\e\|_{L^2}),
\ee
where $\tilde{\rho}=\rho_1+i\rho_2$, and where we used $|\eta|+b^2\lesssim \l^2$ from \eqref{dataatteta0} \eqref{aprioribound}.\\

{\bf step 2} Control of $\tilde{R}$, $R(\e)$ and $\Psi_{\eta,b}$.\\

Note first that \eqref{cnofeof} yields:
\be\label{bins1}
\|\e_z\|_{L^2}\lesssim \alpha^*,\,\|\nabla\e_z\|_{L^2}\lesssim \alpha^*\l.
\ee
Furthermore, \fref{pointwiseboudnut} yields:
\be\label{bins2}
\|\e_ze^{-|y|}\|_{L^2}+\|\e_ze^{-|y|}\|_{L^\infty}\lesssim\alpha^* \l^{m+2}.
\ee
Thus, the terms involving $\tilde{R}$ given by \fref{defrtide} are estimated using \eqref{bins1} and \fref{bins2}:
\be
\label{estrtilde}
\|\tilde{R}\|_{L^2}\lesssim \alpha^* \lambda^{m+2}.
\ee 
The terms involving $R(\e)$ are easily estimated using Sobolev, \eqref{bins1}, \eqref{bins2} and the bootstrap bounds \fref{aprioribound}:
\be
\label{estre}
\int |R(\e)|e^{-|y|}\lesssim\|\e_ze^{-c|y|}\|_{L^\infty}\|\e\|_{L^2}+\|\e\|_{L^2}^2+\|\e\|_{H^1}^3\lesssim K\l^{m+3}.
\ee
Eventually, from \fref{mercredi2} and \fref{dataatteta0}:
\be
\label{estppsietab}
\|\Psi_{\eta,b}\|_{L^2}\lesssim e^{-\frac{c}{\sqrt{\eta}}}\lesssim \alpha^*\l^{m+2},
\ee
for $|\eta|<\eta(\alpha^*)$ small enough.\\

{\bf step 3} The law of $b$.\\

We take the inner product of the equation \eqref{eqeonereal} of $\e_1$ with $-\Lambda\T$ and we sum it with the inner product of equation \eqref{eqeoneim} of $\e_2$ with $\Lambda\S$.   We obtain after integrating by parts:
\be\label{claw1}
|b_s+b^2+\eta|\lesssim |\Re(\e,\overline{Q_{\eta,b}})| +\l^2\|\e\|_{L^2}+|Mod(t)|\|\e\|_{L^2}+(K\l +\alpha^*)\l^{m+2}
\ee
where we used the second orthogonality condition in \eqref{orthee}, the computation of the inner product \eqref{pscal1}, and \eqref{estrtilde}-\eqref{estppsietab}.\\

{\bf step 4} The law of $\l$.\\

We take the inner product of the equation \eqref{eqeonereal} of $\e_1$ with $|y|^2\S$ and we sum it with the inner product of equation \eqref{eqeoneim} of $\e_2$ with $|y|^2\T$.   We obtain after integrating by parts:
\be
\label{claw2}
\left|\lsl+b\right|\lesssim \l^2\|\e\|_{L^2}+|Mod(t)|\|\e\|_{L^2}+(K\l+\alpha^*)\l^{m+2},
\ee
where we used the first orthogonality condition in \eqref{orthee}, the computation of the inner product \eqref{pscal2}, and \eqref{estrtilde}-\eqref{estppsietab}.\\

{\bf step 5} The law of $\tgamma$.\\

We take the inner product of the equation \eqref{eqeonereal} of $\e_1$ with $\rho_2$ and we sum it with the inner product of equation \eqref{eqeoneim} of $\e_2$ with $-\rho_1$.   We obtain after integrating by parts:
\be\label{claw3}
|\tgamma_s|\lesssim |b_s+b^2+\eta|+\l^2\|\e\|_{L^2}+|Mod(t)|\|\e\|_{L^2}+(K\l+\alpha^*)\l^{m+2}
\ee
where we used the third orthogonality condition in \eqref{orthee}, the computation of the inner product \eqref{pscal3}, and \eqref{estrtilde}-\eqref{estppsietab}.

In view of \eqref{claw1}-\eqref{claw3}, we obtain:
\be
\label{paprmatersmeiux}
\left|\lsl+b\right|\lesssim \l^2\|\e\|_{L^2}+(K\l+\alpha^*)\l^{m+2},
\ee
and
\be
\label{paprmatersmeiuxbis}
 |b_s+b^2+\eta|+|\tgamma_s|\lesssim |\Re(\e,\overline{Q_{\eta,b}})|+\l^2\|\e\|_{L^2}+(K\l+\alpha^*)\l^{m+2}.
\ee
Furthermore, we have:
\bea
\label{domani1} |\Re(\e,\overline{Q_{\eta,b}})|&\lesssim& |(\e_1,Q)|+(|\eta|+|b|)\|\e\|_{L^2}\\
\nonumber &\lesssim& |(\e_1,Q)|+K\l^{m+2},
\eea
where we used \eqref{uniformdeacy}, \eqref{defqetab}, and the fact that $|\eta|+|b|\lesssim \l$ from \eqref{aprioribound}, \eqref{parametersinitiallyone}. The estimates \fref{controlmodulation} and \fref{controlmodulationbis}  then follow from the bootstrap bound \eqref{aprioribound}, \fref{paprmatersmeiux}, \fref{paprmatersmeiuxbis}, \eqref{domani1}, and the degeneracy \fref{degeneracyeunqlemma}.\\

{\bf step 6} Proof of the degeneracy \fref{degeneracyeunqlemma}.\\

We take the inner product of the equation \eqref{eqeonereal} of $\e_1$ with $\S$ and we sum it with the inner product of equation \eqref{eqeoneim} of $\e_2$ with $\T$.   We obtain after integrating by parts:
\bea
\label{estderiveeeq}
 \left|\frac{d}{ds}\left(Re(\e,\overline{Q_{\eta,b}})\right)\right| &  \lesssim  & \l^2\|\e\|_{L^2}+|Mod(t)|\|\e\|_{L^2}+(K\l+\alpha^*)\l^{m+2}\\
 & \lesssim & (K\l+\alpha^*)\l^{m+2}
 \eea
where we used the computation of the inner product \eqref{pscal0}, \eqref{estrtilde}-\eqref{estppsietab},  and the bootstrap bound \fref{aprioribound}. Now from \fref{aprioribound}, \fref{parametersinitiallyone}:
$$|t|+|b|\lesssim \lambda(t),$$ which yields: $\forall p\geq 0$, 
\be
\label{cnoefeoif}
\int_{t}^{T_{\eta}}(\lambda(\tau))^pd\tau\lesssim  \int_{t}^{T_{\eta}}(\lambda_{\eta}(\tau))^pd\tau\lesssim |t|(\lambda_{\eta}(t))^p\lesssim (\lambda(t))^{p+1}.
\ee
Thus the time integration from $t$ to $T_{\eta}$ of \fref{estderiveeeq} yields:
 $$|Re(\e,\overline{Q_{\eta,b}})(t)|\lesssim \int_t^{T_{\eta}}\left[K(\lambda(\tau))^{m+1}+\alpha^*(\l(\tau))^{m}\right]d\tau\lesssim K\l(t)^{m+2}+\alpha^*\l(t)^{m+1},$$
 which together with \eqref{domani1} concludes the proof of \fref{degeneracyeunqlemma} and of Lemma \ref{lemmapreliminary}.

%%%%%%%%%%%%%%%%%%%%%%%%%%%%%%%%% 
%%%%%%%%%%%%%%%%%%%%%%%%%%%%%%%%%

\subsection{Mixed energy/Morawetz Lyapounov functional}

%%%%%%%%%%%%%%%%%%%%%%%%%%%%%%%%%
%%%%%%%%%%%%%%%%%%%%%%%%%%%%%%%%%

We first rewrite the a priori bound \fref{aprioribound} as follows:
\be
\label{aprioritu}
\|\nabla \ut\|_{L^2}\leq K \lambda^m, \ \ \|\ut\|_{L^2}\leq K\lambda^{m+1}.
\ee
We will also use the a priori bounds from \fref{aprioribound}, \fref{parametersinitiallyone} and \eqref{dataatteta0}:
\be
\label{aprioriboundblambda}
\eta+b^2\lesssim \l^2.
\ee

We let $A>0$ be a large enough constant which will be chosen later and let $\phi:\RR^2\goto \RR$ be a smooth radially symmetric cut off function with 
\be
\label{defphi}
\phi'(r)=\left\{\begin{array}{ll} r \ \ \mbox{for} \ \ r\leq 1,\\ 3-e^{-r} \ \ \mbox{for}\ \ r\geq 2.\end{array}\right .
\ee
Let $$F(u)=\frac{1}{4}|u|^4, \ \ f(u)=u|u|^2 \ \ \mbox{so that} \ \ F'(u)\cdot h=Re(f(u)\overline{h}).$$ We claim the following generalized energy estimate on the linearized flow \fref{un3}:

\begin{Lemma}[Algebraic generalized energy/Morawetz estimate]
\label{lemma:timederivative}
Let $u=w+\tu$ where $w$ satisfies the bound \eqref{aprioriboundwsw}, and $\tu$ satisfies \eqref{un3} and the a priori bound 
\be
\label{aprioritubis}
\|\nabla \ut\|_{L^2}\lesssim \lambda, \ \ \|\ut\|_{L^2}\lesssim \l^2.
\ee
Let $b, \l$ satisfying the bounds \eqref{aprioriboundblambda} and
\be
\label{roughboundmodulationbis}
|Mod(t)|\lesssim \l^2.
\ee
Let
\bea
\label{defI}
\nonumber \mathcal I(\ut) & = & \frac{1}{2}\int |\n\tu|^2 +\half\int \frac{|\tu|^2}{\l^2}-\int \left[F(w+\tu)-F(w)-F'(w)\cdot\tilde{u}\right]\\
 & + &  \half\frac{b}{\l}\Im\left(\int A\nabla\phi\left(\frac{x}{A\l}\right)\cdot\nabla\tu\overline{\tu}\right),
\eea
\bea
\label{kutilde}
 \mathcal J(\ut) & = & -\frac{1}{\l^2}\Im\left(\int w^2\overline{\tu}^2\right)-\Re\left(\int w_t\overline{(2|\tu|^2w+\tu^2\overline{w})}\right)\\
\nonumber &+ &\frac{b}{\l^2}\left\{\int \frac{|\tu|^2}{\l^2}+\Re\left(\int\nabla^2\phi\left(\frac{x}{A\l}\right)(\n\tu,\overline{\n\tu})\right)-\frac{1}{4A^2}\left(\int\Delta^2\phi\left(\frac{x}{A\l}\right)\frac{|\tu|^2}{\l^2}\right)\right\}\\
\nonumber & + & \frac{b}{\l}\Re\left(\int A\nabla\phi\left(\frac{x}{A\l}\right)(2|\tu|^2w+\tu^2\overline{w})\cdot\overline{\nabla w}\right),
\eea
\bea
\label{limearterm}
\nonumber \mathcal L(\ut) & = &  \Im\left\{\int\left[\Delta\psi-\frac{\psi}{\l^2}+(2|w|^2\psi-w^2\overline{\psi})+i\frac{b}{\l}A\nabla\phi\left(\frac{x}{A\l}\right)\cdot\nabla\psi\right.\right.\\
&+ & \left.\left. i\frac{b}{2\l^2}\Delta\phi\left(\frac{x}{A\l}\right)\psi\right]\overline{\tu}\right\},
\eea
then there holds:
\bea
\label{crc6}
\frac{d}{dt}\mathcal I(\ut) = \mathcal J(\ut)+\mathcal L(\ut)+ O\left(\l^2\|\psi\|_{L^2}^2+\frac{\|\tu\|_{L^2}^2}{\l^2}+\|\tu\|_{H^1}^2\right).
\eea
\end{Lemma}

\begin{remark} The virtue of \fref{crc6} is to keep track of the {\em quadratic terms} $\mathcal I(\ut), \mathcal J(\ut)$ which will turn out to involve to leading order the {\it coercive} quadratic form of the linearized energy near $Q$ given by $(L\ut,\ut)$, and of the leading order {\em linear} term $\mathcal L(\ut)$ involving the error $\psi$.
\end{remark}

{\bf Proof of Lemma \ref{lemma:timederivative}} This Lemma is very similar to Lemma 3.3 in \cite{RS} and we briefly recall the proof for the reader's convenience.\\

{\bf step 1} Algebraic derivation of the energetic part.\\

We compute from \fref{un3}:
\bea
\label{cpepjpvudpov}
& & \frac{d}{dt}\bigg\{\frac{1}{2}\int |\n\tu|^2 +\half\int \frac{|\tu|^2}{\l^2}-\int \left[(F(u)-F(w)-F'(w)\cdot\tu)\right]\bigg\}\\
\nonumber & = &- \Re\left(\partial_t\tu,\overline{\Delta \tu-\frac{1}{\lambda^2}\tu+(f(u)-f(w))}\right)-\frac{\lambda_t}{\lambda^3}\int|\tu|^2\\
\nonumber & - & \Re\left(\partial_tw,(\overline{f(\tu+w)-f(w)-f'(w)\cdot\tu)}\right)\\
\nonumber & = & \Im\left(\psi,\overline{\Delta \tu-\frac{1}{\lambda^2}\tu+(f(u)-f(w))}\right)-\frac{1}{\lambda^2}\Im\left((f(u)-f(w)),\overline{\ut}\right)\\
\nonumber& - & \frac{\lambda_t}{\lambda^3}\int|\tu|^2-\Re\left(\partial_tw,\overline{(f(\tu+w)-f(w)-f'(w)\cdot\tu)}\right)\\
\nonumber& = & \Im\left(\psi,\overline{\Delta \tu-\frac{1}{\lambda^2}\tu+(2|w|^2\ut+\overline{\ut}w^2)}\right)-\frac{1}{\lambda^2}\Im\int \overline{\ut}^2w^2\\
\nonumber& - & \frac{\lambda_t}{\lambda^3}\int|\tu|^2-\Re\left(\partial_tw,\overline{(\overline{w}\tu^2+2w|\tu|^2)}\right)\\
\nonumber& + & \Im\left(\psi-\frac{1}{\lambda^2}\ut,\overline{(f(w+\ut)-f(w)-f'(w)\cdot\ut)}\right)-\Re\left(\partial_tw,\overline{\ut|\ut|^2}\right)
\eea
where we used that $f'(w)\cdot\ut=2|w|^2\ut+w^2\overline{\ut}.$ We first estimate from \fref{roughboundmodulationbis}:
\be
\label{neiohoeghe}
-\frac{\lambda_t}{\lambda^3}\int|\tu|^2=\frac{b}{\lambda^4}\int|\ut|^2-\frac{1}{\lambda^4}\left(\lsl+b\right)\|\ut\|_{L^2}^2=\frac{b}{\lambda^4}\int|\ut|^2+O\left(\frac{\|\tu\|_{L^2}^2}{\l^2}\right).
\ee
It remains to estimate the last line in the RHS \fref{cpepjpvudpov}. For the quadratic and higher terms, we estimate using the a priori bounds \fref{aprioriboundwsw}, \fref{aprioritubis}:
\bea
\label{crc3}
\nonumber & & \left|\Im\left(\psi-\frac{1}{\lambda^2}\ut,\overline{(f(w+\ut)-f(w)-f'(w)\cdot\ut)}\right)\right|\\
\nonumber& = & \left|\Im\left(\psi-\frac{1}{\lambda^2}\ut,\overline{(\ut^2\overline{w}+2|\tu|^2w+|\tu|^2\tu)}\right)\right|\\
\nonumber& \lesssim & \|\psi\|_{L^2(\RR^2)}\|\tu\|_{L^6}^2(\|w\|_{L^6}+\|\tu\|_{L^6})+\frac{(1+\|w\|_{L^2})}{\l^2}\|\ut\|_{L^6}^3\\
 & \lesssim& \l^2\|\psi\|_{L^2}^2+\|\tu\|_{H^1}^2.
\eea
For the cubic term hitting $w_t$, we replace $w_t$ using \eqref{eqwgobale}, integrate by parts and use \fref{aprioriboundwsw}, \fref{aprioritubis} to estimate:
\bea
\label{cnkeneofho}
\nonumber \left|\int w_t\overline{|\tu|^2\tu}\right|& \lesssim & \|w\|_{H^{3/2}}\||\tu|^2\tu\|_{H^{1/2}(\RR^2)}+\|w\|^3_{L^6}\|\tu\|^3_{L^6}
+\|\psi\|_{L^2(\RR^2)}\|\tu\|^3_{L^6}\\
\nonumber & \lesssim & \frac{1}{\l^{3/2}}\|\tu\|^{1/2}_{L^2}\|\tu\|^{5/2}_{H^1}+\frac{1}{\lambda^2}\|\ut\|_{H^1}^2\|\tu\|_{L^2}+\|\psi\|_{L^2}\|\tu\|^{2}_{H^1}\|\tu\|_{L^2}\\
& \lesssim & \l^2\|\psi\|_{L^2}^2+\|\tu\|_{H^1}^2.
\eea
Injecting \fref{neiohoeghe}, \fref{crc3}, \fref{cnkeneofho} into \fref{cpepjpvudpov} yields the preliminary computation:
\bea
\label{cnkoheofh}
 & &  \frac{d}{dt}\bigg\{\frac{1}{2}\int |\n\tu|^2+\half\int \frac{|\tu|^2}{\l^2}- \int (F(u)-F(w)-F'(w)\cdot\tu)\bigg\}\\
\nonumber & = & -\frac{1}{\l^2}\Im\left(\int w^2\overline{\tu}^2\right)-\Re\left(\int w_t\overline{(2|\tu|^2w+\tu^2\overline{w})}\right)+\frac{b}{\l^2}\int \frac{|\tu|^2}{\l^2}\\
\nonumber & + & \Im\left(\int\left[\Delta\psi-\frac{\psi}{\l^2}+(2|w|^2\psi-w^2\overline{\psi})\right]\overline{\tu}\right)\\
\nonumber & + & O\left(\l^2\|\psi\|_{L^2}^2+\frac{\|\tu\|_{L^2}^2}{\l^2}+\|\tu\|_{H^1}^2\right).
\eea

{\bf step 2} Algebraic derivation of the localized virial part.\\

Let $$\nabla\tilde{\phi}(t,x)=\frac{b}{\l}A\nabla\phi\left(\frac{x}{A\l}\right),$$ then
\bea
\label{firsttermvirieloc}
 & & \half\frac{d}{dt}\left\{\frac{b}{\l}\Im\left(\int A\nabla\phi\left(\frac{x}{A\l}\right)\nabla\tu\overline{\tu}\right)\right\} \\
 \nonumber & =&  \half \Im\left(\int\partial_t\nabla \tilde{\phi}\cdot \nabla\ut\overline{\ut}\right)+  \Re\left(\int i\partial_t\ut\left[\overline{\half \Delta\tilde{\phi}\tu+\nabla\tilde{\phi}\cdot\nabla\tilde{u}}\right]\right).
\eea
We estimate in brute force using \fref{defphi}, \fref{roughboundmodulationbis}, \fref{aprioriboundblambda}:
$$
\left|\partial_t\nabla \tilde{\phi}\right|\lesssim \frac{A}{\lambda^3}\left(|b_s|+b^2+b\left|\lsl+b\right|\right)\lesssim \frac{A}{\l}
$$
from which: 
\be
\label{calculprline}
\left|\Im\left(\int\partial_t\nabla \tilde{\phi}\cdot \nabla\ut\overline{\ut}\right)\right|\lesssim \frac{A}{\lambda}\|\ut\|_{L^2}\|\nabla\ut\|_{L^2}\lesssim \frac{\|\tu\|_{L^2}^2}{\l^2}+\|\tu\|_{H^1}^2.
\ee
The second term in \fref{firsttermvirieloc} corresponds to the localized Morawetz multiplier, and we get from \fref{un3} and the classical Pohozaev integration by parts formula:
\bea
\label{nceohoeoud}
& & \Re\left(\int i\partial_t\ut\left[\overline{\half \Delta\tilde{\phi}\tu+\nabla\tilde{\phi}\cdot\nabla\tilde{u}}\right]\right)=\\
\nonumber & = &   \frac{b}{\l^2}\Re\left(\int\nabla^2\phi\left(\frac{x}{A\l}\right)(\n\tu,\overline{\n\tu})\right)-\frac{1}{4}\frac{b}{A^2\l^4}\left(\int\Delta^2\phi\left(\frac{x-\a}{A\l}\right)|\tu|^2\right)\\
\nonumber &  - & \frac{b}{\l}\Re\left(\int A\nabla\phi\left(\frac{x}{A\l}\right)(|u|^2u-|w|^2w)\cdot\overline{\nabla\tu}\right) -  \half\frac{b}{\l^2}\Re\left(\int\Delta\phi\left(\frac{x}{A\l}\right)(|u|^2u-|w|^2w)\overline{\tu}\right)\\
\nonumber & - &\frac{b}{\l}\Re\left(\int A\nabla\phi\left(\frac{x}{A\l}\right)\psi\cdot\overline{\nabla\tu}\right)\ds -\half\frac{b}{\l^2}\Re\left(\int\Delta\phi\left(\frac{x}{A\l}\right)\psi\overline{\tu}\right).
\eea
We now expand the nonlinear terms and estimate the cubic and higher terms:
\bea
\label{ceheoe}
\nonumber & &  \left|-\frac{b}{\l}\Re\left(\int A\nabla\phi\left(\frac{x}{A\l}\right)(2|\tu|^2w+\tu^2\overline{w}+|\tu|^2\tu)\cdot\overline{\nabla\tu}\right)\right .\\
\nonumber & - & \left .\half\frac{b}{\l^2}\Re\left(\int\Delta\phi\left(\frac{x}{A\l}\right)(2|\tu|^2w+\tu^2\overline{w}+|\tu|^2\tu)\overline{\tu}\right)\right|\\
\nonumber & \lesssim &(\|\tu\|^3_{L^6(\RR^2)}+\|\tu\|^2_{L^6(\RR^2)}\|w\|_{L^6(\RR^2)})\|\n\tu\|_{L^2(\RR^2)}+\frac{1}{\l}(\|\tu\|^4_{L^4(\RR^2)}+\|\tu\|^3_{L^4(\RR^2)}\|w\|_{L^4(\RR^2)})\\
& \lesssim & \frac{\|\tu\|_{L^2}^2}{\l^2}+\|\tu\|_{H^1}^2
\eea
where we have used \eqref{aprioritubis} and \eqref{aprioriboundwsw}. The remaining quadratic terms in \fref{nceohoeoud} are integrated by parts:
\bea
\label{lpusfpfwu}
\nonumber &- &\frac{b}{\l}\Re\left(\int A\nabla\phi\left(\frac{x}{A\l}\right)\psi\cdot\overline{\nabla\tu}\right)\ds -\half\frac{b}{\l^2}\Re\left(\int\Delta\phi\left(\frac{x}{A\l}\right)\psi\overline{\tu}\right)\\
& = & \Im \left(\int \left[i\frac{b}{\lambda}A\nabla\phi\left(\frac{x}{A\l}\right)\cdot\nabla \psi+i\frac{b}{2\l^2}\Delta\phi\left(\frac{x-\a}{A\l}\right)\psi\right]\overline{\ut}\right),
\eea
and
\bea
\label{bingo}
\nonumber & - & \frac{b}{\l}\Re\left(\int A\nabla\phi\left(\frac{x}{A\l}\right)(2|w|^2\ut+w^2\overline{\ut})\cdot\overline{\nabla\tu}\right)\\
\nonumber & - &  \half\frac{b}{\l^2}\Re\left(\int\Delta\phi\left(\frac{x}{A\l}\right)(2|w|^2\ut+w^2\overline{\ut})\overline{\tu}\right)\\
& = &  \frac{b}{\l}\Re\left(\int A\nabla\phi\left(\frac{x}{A\l}\right)(2|\tu|^2w+\tu^2\overline{w})\cdot\overline{\nabla w}\right).
\eea
Injecting \eqref{ceheoe}, \fref{lpusfpfwu}, \eqref{bingo} into \fref{nceohoeoud} yields after a further integration by parts:
\bee
& & \Re\left(\int i\partial_t\ut\left[\overline{\half \Delta\tilde{\phi}\tu+\nabla\tilde{\phi}\cdot\nabla\tilde{u}}\right]\right)\\
\nonumber & = &   \frac{b}{\l^2}\Re\left(\int\nabla^2\phi\left(\frac{x}{A\l}\right)(\n\tu,\overline{\n\tu})\right)-\frac{1}{4}\frac{b}{A^2\l^4}\left(\int\Delta^2\phi\left(\frac{x}{A\l}\right)|\tu|^2\right)\\
\nonumber & + & \frac{b}{\l}\Re\left(\int A\nabla\phi\left(\frac{x}{A\l}\right)(2|\tu|^2w+\tu^2\overline{w})\cdot\overline{\nabla w}\right)\\
\nonumber & + & \Im\left(\int\left[i\frac{b}{\l}A\nabla\phi\left(\frac{x}{A\l}\right)\cdot\nabla\psi+i\frac{b}{2\l^2}\Delta\phi\left(\frac{x}{A\l}\right)\psi\right]\overline{\tu}\right)+O\left(\frac{\|\tu\|_{L^2}^2}{\l^2}+\|\tu\|_{H^1}^2\right).
\eee
We now inject this together with \fref{calculprline} into \fref{firsttermvirieloc} which together with \fref{cnkoheofh} concludes the proof of \fref
{crc6}.

%%%%%%%%%%%%%%%%%%%%%%%%%%%%%%%%%%%%%%%%%%%%
%%%%%%%%%%%%%%%%%%%%%%%%%%%%%%%%%%%%%%%%%%%%

\subsection{Proof of the bootstrap Lemma \ref{bootlemma}}

%%%%%%%%%%%%%%%%%%%%%%%%%%%%%%%%%%%%%%%%%%%%
%%%%%%%%%%%%%%%%%%%%%%%%%%%%%%%%%%%%%%%%%%%%

We are now in position to close the bootstrap estimates \fref{estktun} as a consequence of the Lyapounov control of Lemma \ref{lemma:timederivative}.\\

{\bf step 1}  Coercitivity of $\mathcal I$.\\

We first claim the coercitivity of $\mathcal I(\ut)$ given by \fref{defI}: 
\be
\label{coeritilde}
\mathcal I(\ut)\geq c_0\left(\|\nabla \ut\|_{L^2}^2+\frac{1}{\l^2}\|\ut\|_{L^2}^2\right)+O\left((K^2\l^2+\alpha^*)\l^{2m}\right)
\ee
for some universal constant $c_0>0$. Indeed, we first renormalize:
\bee
\nonumber \mathcal I(\ut)& = & \frac{1}{2\l^2}\left\{\int |\n\e|^2 +\int |\e|^2-2\int \left(F(Q_{\eta,b}+\e_z+\e)-F(Q_{\eta,b}+\e_z)-F'(Q_{\eta,b}+\e_z)\cdot\e\right)\right .\\
& + &\left . b\Im\left(\int A\nabla\phi\left(\frac{y}{A}\right)\nabla\e\overline{\e}\right)\right\}.
\eee
We then estimate by homogeneity:
\bee
%\label{esthomogeniety}
& & 2\int\left[F(Q_{\eta,b}+\e_z+\e)-F(Q_{\eta,b}+\e_z)-F'(Q_{\eta,b}+\e_z)\cdot \e\right] \\
&  =  &\int\left[ Re(\e^2\overline{W}^2)+2|\e|^2|W|^2\right]+ O\left(\int |W||\e|^3+|\e|^4\right)\\
& = & 3\int\e_1^2Q^2+\int\e_2^2Q^2+O\left((|\eta|+|b|+\|\e_ze^{-c|y|}\|_{L^2}+\|\e_z\|_{L^4}^2)\|\e\|_{H^1}^2+\|\e\|_{H^1}^3+\|\e\|_{H^1}^4\right)\\
& = & 3\int\e_1^2Q^2+\int\e_2^2Q^2+O((\l+\a^*)\|\e\|_{H^1}^2)
\eee
where we used the proximity of $Q_{b,\eta}$ to $Q$ given by \fref{uniformdeacy} \eqref{defqetab}, the bound \eqref{aprioriboundblambda} for $\eta$ and $b$, the a priori bound \eqref{aprioribound} for $\e$, and the estimates \eqref{bins1} \eqref{bins2} for $\e_z$. We thus obtain:
\be
\label{esticoineo}
\mathcal I(\ut)\geq \frac{1}{\l^2}\left\{(L_+\e_1,\e_1)+(L_-\e_2,\e_2)+O((\l+\alpha^*)\|\e\|_{H^1}^2)\right\}.
\ee
 We now recall the following coercivity property of the linearized energy which is a well known consequence of the variational characterization of $Q$:

\begin{Lemma}[Coercivity of the linearized energy, \cite{W2}, \cite{MR1}, \cite{MR4}]
\label{lemmacoerc}
There holds for some universal constant $c_0>0$ :  $\forall \e\in H^1$ radially symmetric, 
\bea
\label{coerclinearenergy}
& & (L_+\e_1,\e_1)+(L_-\e_2,\e_2) \geq c_0  \|\e\|_{H^1}^2\\
\nonumber & - & \frac{1}{c_0}\left\{(\e_1,Q)^2+(\e_1,|y|^2Q)^2+(\e_2,\rho)^2\right\}. \eea
\end{Lemma}

The a priori bound \eqref{aprioribound} and the degeneracy \fref{degeneracyeunqlemma} imply:
\be
\label{newdegeneracy}
(\e_1,Q)^2\lesssim\left((K\l+\alpha^*)\l^{m+1}\right)^2\lesssim (K^2\l^2+\alpha^*)\lambda^{2m+2}.
\ee
Injecting the choice of orthogonality conditions \fref{orthee} and \eqref{newdegeneracy} into \fref{coerclinearenergy} and then into \fref{esticoineo} now yields \fref{coeritilde} for $\alpha^*$ small enough.\\

{\bf step 2} Coercitivity of $\mathcal J(\ut)$\\

Let $\mathcal J(\ut)$ be given by \fref{kutilde}, we claim:
\be
\label{firtboundq}
\mathcal{J}(\tu)\geq c\frac{b}{\lambda^4}\left(\int|\nabla \e|^2e^{-\frac{|y|}{\sqrt{A}}}+\int|\e|^2\right)+O\left((K^2\l+\alpha^*)\lambda^{2m-1}\right)
\ee
for some universal constant $c>0$. Recall the decomposition \fref{defqtilde}, we first claim:
\be
\label{replacementkone}
\mathcal{J}(\ut)=\mathcal{J}_1(\ut)+O\left(A\frac{\|\tu\|^2_{L^2}}{\l^2}+\|\tu\|^2_{H^1}\right)
\ee with 
\bea
\label{defkoneutilde}
\mathcal J_1(\ut) & = & -\frac{1}{\l^2}\Im\left(\int \tq^2\overline{\tu}^2\right)-\Re\left(\int \tq_t\overline{(2|\tu|^2\tq+\tu^2\overline{\tq})}\right)\\
\nonumber &+ &\frac{b}{\l^2}\left\{\int \frac{|\tu|^2}{\l^2}+\Re\left(\int\nabla^2\phi\left(\frac{x}{A\l}\right)(\n\tu,\overline{\n\tu})\right)-\frac{1}{4A^2}\left(\int\Delta^2\phi\left(\frac{x}{A\l}\right)\frac{|\tu|^2}{\l^2}\right)\right .\\
\nonumber & + & \left .\l\Re\left(\int A\nabla\phi\left(\frac{x}{A\l}\right)(2|\tu|^2\tq+\tu^2\overline{\tq})\cdot\overline{\nabla \tq}\right)\right\}.
\eea
Indeed, we expand $w=\tq+z$ and estimate the remaining terms using \fref{cnofeof}, \eqref{pointwiseboudnut}, \eqref{aprioriboundwsw} and \eqref{bins2}. We obtain:
$$\frac{1}{\l^2}\int(|\tq||z|+|z|^2)|\ut|^2\lesssim \frac{1}{\l^2}(\|\e_ze^{-c|y|}\|_{L^2}\|\ut\|_{L^4}^2+\|z\|_{L^\infty}^2\|\ut\|^2_{L^2})\lesssim \frac{\|\tu\|^2_{L^2}}{\l^2}+\|\tu\|^2_{H^1},
$$
\bee
&&\left|\int \pa_tz(\overline{2|\ut|^2w+\ut^2\overline{w}})\right|+\left|\int \pa_t\tq(\overline{2|\ut|^2z+\ut^2\overline{z}})\right|\\
&\lesssim& \|\pa_tz\|_{L^4}\|\ut\|^2_{L^4}\|w\|_{L^4}+\frac{1}{\l^2}\|\e_ze^{-c|y|}\|_{L^2}\|\ut\|_{L^4}^2 \lesssim \frac{\|\tu\|^2_{L^2}}{\l^2}+\|\tu\|^2_{H^1},
\eee
$$\frac{b^2}{\l^3}A\int |\ut|^2|z||\nabla z|\lesssim \frac{A}{\l}\|\ut\|_{L^2}^2\lesssim \frac{\|\tu\|^2_{L^2}}{\l^2},$$
and
\bee
&& \left|\frac{b}{\l}\Re\left(\int A\nabla\phi\left(\frac{x}{A\l}\right)\left[(2|\tu|^2\tq+\tu^2\overline{\tq})\cdot\overline{\nabla z}+(2|\tu|^2z+\tu^2\overline{z})\cdot\overline{\nabla w}\right]\right)\right|\\
&\lesssim&\frac{A}{\l}\|\ut\|^2_{L^2}\left(\frac{1}{\l}[1+\|\nabla z\|^2_{L^\infty}+\|z\|^2_{L^{\infty}}]+\frac{1}{\lambda^3}\|\e_ze^{-c|y|}\|_{L^\infty}\right)\lesssim A\frac{\|\tu\|^2_{L^2}}{\l^2},
\eee
and \fref{replacementkone} follows.\\
We now expand the $\pa_t\tq$ term from \fref{defqtilde}, \fref{roughboundmodulationbis}, \eqref{aprioriboundblambda}:
\bee
\tq_t & = & \left(\frac{i}{\l^2}+\frac{b}{\l^2}\right)\tq+\frac{b}{\l}\frac{x}{\lambda}\cdot\nabla \tq+O\left(\frac{|Mod(t)|+b^2+\eta}{\l^3}e^{-c\frac{|x|}{\l}}\right)\\
& = & \left(\frac{i}{\l^2}+\frac{b}{\l^2}\right)\tq+\frac{b}{\l}\frac{x}{\lambda}\cdot\nabla \tq+O\left(\frac{1}{\l}\right),
\eee
which yields:
\bee
\nonumber &- & \Re\left(\int \tq_t\overline{(2|\tu|^2\tq+\tu^2\overline{\tq})}\right)\\
&=& \ds\frac{1}{\l^2}\Im\left(\int\tq\overline{(2|\tu|^2\tq+\tu^2\overline{\tq})}\right)- \frac{b}{\l^2}\Re\left(\int (2|\tu|^2\tq+\tu^2\overline{\tq})\overline{\tq}\right)\\
&& -\frac{b}{\l}\Re\left(\int \frac{x}{\l}(2|\ut|^2\tq+\ut^2\overline{\tq})\cdot\overline{\nabla\tq}\right)+  O\left(\frac{1}{\l}\|\tq\|_{L^\infty}\|\ut\|^2_{L^2}\right)\\
&=& \ds\frac{1}{\l^2}\Im\left(\int\tq\overline{(2|\tu|^2\tq+\tu^2\overline{\tq})}\right)- \frac{b}{\l^2}\Re\left(\int (2|\tu|^2\tq+\tu^2\overline{\tq})\overline{\tq}\right)\\
&& -\frac{b}{\l}\Re\left(\int \frac{x}{\l}(2|\ut|^2\tq+\ut^2\overline{\tq})\cdot\overline{\nabla\tq}\right)+  O\left(\frac{\|\tu\|^2_{L^2}}{\l^2}\right).
\eee
We inject this estimate into \fref{defkoneutilde} and write the result in renormalized variables:
\bea
\label{sunday}\mathcal{J}_1(\tu)
&= & \frac{b}{\l^4}\left\{\Re\left(\int\nabla^2\phi\left(\frac{y}{A}\right)(\n\e,\overline{\n\e})\right)+\int|\e|^2 \right .\\
\nonumber&- &\left . \int \left[(3\S^2+\T^2)\e_1^2+4\S\T\e_1\e_2+(\S^2+3\T^2)\e_2^2\right]- \frac{1}{4A^2}\int\Delta^2\phi\left(\frac{y}{A}\right)|\e|^2\right .\\
\nonumber& + & \left.\Re\left(\int \left(A\nabla\phi\left(\frac{y}{A}\right)-y\right)(2|\e|^2Q_{\eta,b}+\e^2\overline{Q_{\eta,b}})\cdot\overline{\nabla Q_{\eta,b}})\right)\right\}+   O\left(\frac{\|\tu\|^2_{L^2}}{\l^2}\right).
\eea
From the proximity of $Q_{\eta,b}$ to $Q$ and the control of the {\it full} $L^2$ norm, the above quadratic form is for $A$ large enough a small deformation of the localized in $A$ linearized energy, and hence \fref{coerclinearenergy}, \fref{newdegeneracy} and our choice of orthogonality conditions ensure for $A$ large enough:
\bee
\mathcal{J}(\tu) & \geq & \frac{c_0}{2}\frac{b}{\lambda^4}\left[\int|\nabla \e|^2e^{-\frac{|y|}{\sqrt{A}}}+\int|\e|^2\right]+O\left((K^2\l^2+\alpha^*)\lambda^{2m-1}+A\frac{\|\tu\|^2_{L^2}}{\l^2}+\|\tu\|^2_{H^1}\right)\\
& \geq & O\left((K^2\l+\alpha^*)\lambda^{2m-1}\right)
\eee
where we used \fref{aprioribound} in the last step. This concludes the proof of \fref{firtboundq}.\\

{\bf step 3} Control of $\mathcal L(\ut)$.\\

We now turn to the control of the leading order linear term $\mathcal L(\ut)$ given by \fref{limearterm} and we claim:
\be
\label{estlutilde}
|\mathcal L(\ut)|\lesssim (\l K^2+\alpha^*K)\lambda^{2m-1}.
\ee
We first derive from \fref{defpsiacaculaer}, \eqref{defrtide}, \fref{bins2} and \fref{controlmodulation} the rough bound:
\be
\label{boudpsiretse}
|\Psi|\lesssim |Mod(t)|e^{-c|y|}+|\Psi_{\eta,b}|+|\tilde{R}|\lesssim (K\l+\alpha^*)\l^{m+1}e^{-c|y|}.
\ee
Together with \fref{bins2} and \eqref{aprioribound}, this yields the bound:
\bee
&& \int (|z|^2+|\tq||z|)|\psi||\ut|\lesssim \frac{1}{\l^4}\int (|\e_z|^2+|\e_z|e^{-c|y|})|\Psi||\e|\\
&\lesssim& \frac{1}{\l^4}\|\e_ze^{-c|y|}\|_{L^\infty}(1+\|\e_ze^{-c|y|}\|_{L^\infty})\l^{m+1}\|\e\|_{L^2} \lesssim \alpha^*K\l^{3m}
\eee
and thus $$\mathcal L(\ut)=\mathcal L_1(\ut)+\mathcal L_2(\ut)+O\left((K^2\l+\alpha^*)\lambda^{2m-1}\right)$$ with 
\bea
\label{limeartermone}
\nonumber \mathcal L_j(\ut) & = &  \Im\left\{\int\left[\Delta\psi_j-\frac{\psi_j}{\l^2}+(2|\tq|^2\psi_j-\tq^2\overline{\psi_j})+i\frac{b}{\l}A\nabla\phi\left(\frac{x}{A\l}\right)\cdot\nabla\psi_j\right.\right.\\
&+ & \left.\left. i\frac{b}{2\l^2}\Delta\phi\left(\frac{x}{A\l}\right)\psi_j\right]\overline{\tu}\right\}
\eea
where according to \fref{defpsiacaculaer}: $$\psi_j=\frac{1}{\l^3}\Psi_j\left(\frac{x}{\lambda}\right)e^{i\gamma}, \ \ \Psi_1=\Psi-\tilde{R}+\Psi_{\eta,b}, \ \ \Psi_2=\tilde{R}-\Psi_{\eta,b}.$$
We start with the second term. We estimate in brute force using \fref{estrtilde} \eqref{estppsietab}:
$$\|\psi_2\|_{L^2}\lesssim \frac{1}{\l^2}(\|\Psi_{\eta,b}\|_{L^2}+\|\tilde{R}\|_{L^2})\lesssim\alpha^*\l^m.$$
Now, note that \eqref{pointwiseboudnut} yields the following analog of \eqref{bins2}:
$$\|\nabla\e_ze^{-c|y|}\|_{L^2}\lesssim \alpha^*\l^{m+\frac{5}{2}},$$
which together with \eqref{mercredi2}, \eqref{bins2} yields:
$$\|\nabla\psi_2\|_{L^2}\lesssim \frac{1}{\l^3}(\|\nabla\Psi_{\eta,b}\|_{L^2}+\|\e_ze^{-c|y|}\|_{L^2}+\|\nabla\e_z e^{-c|y|}\|_{L^2})\lesssim \alpha^*\l^{m-1}.$$
These estimates for $\psi_2$ and $\nabla\psi_2$ together with \eqref{aprioriboundblambda} and the a priori bound \eqref{aprioritu} immediately imply:
\bee
|\mathcal L_2(\ut)|&\lesssim& \|\nabla\psi_2\|_{L^2}\|\nabla\ut\|_{L^2}+\frac{1}{\l^2}\|\psi_2\|_{L^2}\|\ut\|_{L^2}+A\|\nabla\psi_2\|_{L^2}\|\ut\|_{L^2}+\frac{1}{\l}\|\psi_2\|_{L^2}\|\ut\|_{L^2}\\
&\lesssim& K\alpha^*\lambda^{2m-1}.
\eee
We now compute $\mathcal L_1(\ut)$ from the explicit formula \fref{defpsiacaculaer} and observe that the corresponding terms lie up to $O(\l^2)$ in the general null space of $L$, and hence a factor $\lambda$ is gained. In other words, we compute after renormalization:
\bea
\label{cneoneogehoeueo}
\nonumber |\mathcal L_1(\ut)| & \lesssim & \frac{|Mod(t)|}{\lambda^4}\left[|(\e_2,L_-(|y|^2Q))|+|(\e_2,L_-Q)|+O(\l\|\e\|_{L^2})\right]\\
\nonumber&+ & \frac{1}{\l^4}\left|\lsl+b\right|\left[|(\e_1,L_+(\Lambda Q))|+O(\l\|\e\|_{L^2})\right]\\
\nonumber & \lesssim & K\frac{\l^{m+2}}{\l^4}\left[K\l+\alpha^*\right]\l^{m+1}+\frac{1}{\l^4}(K+\alpha^*)\l^{m+2}(K\l+\alpha^*)\l^{m+1}\\
& \lesssim & (\l K^2+\alpha^*K)\l^{2m-1}
\eea
where we used the algebraic identities \eqref{structurekernelpair}, the bounds \fref{newdegeneracy}, \fref{controlmodulation}, \eqref{controlmodulationbis}, the orthogonality conditions \eqref{orthee}, and the a priori bound \eqref{aprioribound}. This concludes the proof of \fref{estlutilde}.\\

{\bf step 4} Control of $\ut$.\\

 We derive from \fref{boudpsiretse} the rough bound:
\be
\label{estpsi}
\|\psi\|_{L^2}=\frac{1}{\l^2}\|\Psi\|_{L^2}\lesssim (K\l+\alpha^*)\l^{m-1}.
\ee
We inject the bounds \eqref{aprioritu}, \fref{firtboundq}, \fref{estlutilde} and \eqref{estpsi} into \fref{crc6} and obtain the Lyapounov control:
$$\frac{d}{dt}\mathcal I\gtrsim -(K^2\l+\alpha^*+\alpha^*K)\l^{2m-1}.$$ We integrate this from $t$ to $T_{\eta}$ and use the coercitivity bound \fref{coeritilde} and \fref{cnoefeoif} to conclude:
\bee
\|\nabla \ut(t)\|_{L^2}^2+\frac{\|\ut(t)\|_{L^2}^2}{\l^2}&\lesssim& (K^2\l(t)^2+\alpha^*)\l(t)^{2m}+\int_{t}^{T_{\eta}}[K^2\l(\tau)+\alpha^*+\alpha^*K]\l(\tau)^{2m-1}d\tau\\
&\leq &\left[\frac{K^2}{2}+O(\alpha^*)\right](\lambda(t))^{2m}
\eee
provided $|t_1|$ is sufficiently small. Hence 
\be
\label{improvedboundsquareroot}
K^2\leq \frac{K^2}{2}+O(\alpha^*) \ \ \mbox{and thus} \ \ K\lesssim \sqrt{\alpha^*}.
\ee 

 {\bf step 5} Integration of the law for the parameters.\\
 
 The control of the parameters now follows by reintegrating the modulation equations. Indeed, from \fref{controlmodulation}, \fref{improvedboundsquareroot}:
 \bea
 \label{tobetimeintegrraed}
 \nonumber \left\{\frac{\sqrt{b^2+\eta}}{\l}\right\}_s & = & \frac{b}{\l\sqrt{b^2+\eta}}(b_s+b^2+\eta)-\frac{\sqrt{b^2+\eta}}{\l}\left(\lsl+b\right)\\
 & = &  \frac{b_{\eta}}{\l_{\eta}\sqrt{b_{\eta}^2+\eta}}(b_s+b^2+\eta)+O(\sqrt{\alpha^*}\l_{\eta}^{m+1})=O(\sqrt{\alpha^*}\lambda_{\eta}^m).
 \eea
 We integrate this bound from $t$ to $T_{\eta}$ and recall from \fref{parametersinitiallyone} that $$\frac{\sqrt{b_{\eta}^2+\eta}}{\lambda_{\eta}}=1.$$ Hence: 
 \bea
\nonumber\frac{\sqrt{(b(t))^2+\eta}}{\l(t)}&=&1+O(\sqrt{\alpha^*}\l^m)+O\left(\int_t^{T_{\eta}}\sqrt{\alpha^*}(\lambda_{\eta}(\tau))^{m-2}d\tau\right)\\
\label{reint1}&=& 1+O(\sqrt{\alpha^*}\l_{\eta}^{m-1}).
\eea
\eqref{reint1} implies:
$$b^2+\eta=\l^2+O(\sqrt{\alpha^*}\l_{\eta}^{m+1}),$$
which  together with  \fref{controlmodulation} yields:
$$b_s+\l^2=O(\sqrt{\alpha^*}\l_{\eta}^{m+1})\ \ \mbox{ie} \ \ b_t+1=O(\sqrt{\alpha^*}\l_{\eta}^{m-1}).$$
We integrate from $t$ to $T_{\eta}$ using $b(T_{\eta})+T_{\eta}=b_{\eta}(T_{\eta})+T_{\eta}=0$ and thus:
\be
\label{reint2}
b(t)+t=b(t)-b_{\eta}(t)=\int_{T_\eta}^t\sqrt{\alpha^*}O(\l_{\eta}^{m-1})=O(\sqrt{\alpha^*}\l_{\eta}^m).
\ee
Next, we compare $\l$ with $\l_\eta$. In view of \eqref{reint1}, we have:
\bee
\l&=&\sqrt{b^2+\eta}+O(\sqrt{\alpha^*}\l_{\eta}^m)= \l_\eta+\sqrt{b^2+\eta}-\sqrt{b^2_\eta+\eta}+O(\sqrt{\alpha^*}\l_{\eta}^m)\\
&=& \l_\eta+O(|b-b_\eta|+\sqrt{\alpha^*}\l_{\eta}^m),
\eee
where we used \fref{parametersinitiallyone}. Together with \eqref{reint2}, this yields:
$$\l=\l_\eta+O(\sqrt{\alpha^*}\l_\eta^m) \ \ \mbox{and thus} \ \ \left|\frac{\l}{\l_{\eta}}-1\right|\lesssim \sqrt{\alpha^*}\l^{m-1}.$$
We finally compute the phase using \eqref{controlmodulation}, \eqref{aprioribound} and \eqref{improvedboundsquareroot}:
 $$|\tgamma_s|\lesssim \sqrt{\alpha^*}\l^{m+1} \ \ \mbox{and thus}\ \ \left|\frac{d}{dt}(\gamma-\gamma_{\eta})\right|\lesssim \left|\frac{1}{\l^2}-\frac{1}{\l_{\eta}^2}\right|+\sqrt{\a^*}\l^{m-1}\lesssim \sqrt{\a^*}\l^{m-3},$$
which after integration in time yields:
$$|\gamma-\gamma_{\eta}|\lesssim \sqrt{\alpha^*}\l^{m-2}.$$  Hence, we have obtained the following estimate for the parameters:
\be\label{montag}
\frac{1}{\lambda^{m-1}}\left|\frac{\lambda_\eta}{\lambda}-1\right|+\frac{|b_{\eta}-b|}{\lambda^m}+\frac{|\gamma-\gamma_{\eta}|}{\lambda^{m-2}}\lesssim \sqrt{\alpha^*}.
\ee
 
 {\bf step 6} $H^{\frac 32}$ control.\\
 
 It remains to prove the $H^{\frac 32}$ bound on $[t_1,T_{\eta}]$:
 \be
 \label{htrheehalfbound}
 \|\e\|_{H^{\frac 32}}\lesssim \sqrt{\alpha^*}\l^{m-\frac32}.
 \ee
 In view of \eqref{un3}, 
 $\tilde{u}$ satisfies: 
\be\label{b3demi}
i\p_t\tu+\Delta\tu=-\psi-F_1-F_2,
\ee
where $\psi$ is defined by \eqref{eqwgobale} \eqref{defpsiacaculaer}, $F_1$ is given by:
$$F_1=2|\tq|^2\tu+\tq^2\overline{\tu},$$
$F_2$ is given by:
$$F_2=|u|^2u-|w|^2w-F_1,$$
and where $w$ has been defined in \eqref{defqtilde}. Hence, from standard Strichartz bounds and the smoothing effect of the linear Schr\"odinger flow, there holds:
\be
\label{b3demi1}
\|\nabla^{\frac{3}{2}}\tu\|_{L^\infty_{[t,T_\eta]}L^2}\lesssim
\|\nabla^{\frac{3}{2}}\psi\|_{L^{\frac{4}{3}}_{[t,T_\eta]}L^{\frac{4}{3}}}+\|(1+|x|^2)F_1\|_{L^2_{[t,T_\eta]}H^1}+\|\nabla^{\frac{3}{2}}F_2\|_{L^{\frac{4}{3}}_{[t,T_\eta]}L^{\frac{4}{3}}}.
\ee
In view of the definition of $\psi$ \eqref{eqwgobale} \eqref{defpsiacaculaer} and $\tilde{R}$ \eqref{defrtide}, we have:
\bee
\|\nabla^{\frac{3}{2}}\psi\|_{L^{\frac{4}{3}}}&\lesssim& \frac{1}{\l^3}\left(|Mod(t)|+\|\Psi_{\eta,b}\|_{L^{\frac 43}}+\|\nabla^{\frac{3}{2}}\widetilde{R}\|_{L^{\frac{4}{3}}}\right)\\
&\lesssim& \frac{1}{\l^3}\left(\sqrt{\alpha^*}\l^{m+1}+\|e^{-c|y|}\e_z\|_{L^\infty}+\|e^{-c|y|}\nabla^2\e_z\|_{L^\infty}\right),
\eee
where we used \eqref{mercredi2} \eqref{controlmodulation} \eqref{improvedboundsquareroot} in the last inequality. Now, in view of \eqref{pointwiseboudnut}, we have:
$$\|e^{-c|y|}\e_z\|_{L^\infty}+\|e^{-c|y|}\nabla^2\e_z\|_{L^\infty}\lesssim \alpha^*\l^{m+2},$$ 
which in turn yields:
$$\|\nabla^{\frac{3}{2}}\psi\|_{L^{\frac{4}{3}}}\lesssim\sqrt{\alpha^*} \l^{m-2}.$$
Thus:
\be\label{b3demi2}
\|\nabla^{\frac{3}{2}}\psi\|_{L^{\frac{4}{3}}_{[t,T_\eta]}L^{\frac{4}{3}}}\lesssim\sqrt{\alpha^*} \l^{m-\frac{5}{4}}.
\ee
The $F_1$ term is local in $y$, hence from the a priori bound \eqref{aprioritu} and \eqref{improvedboundsquareroot}, we obtain:
$$\|(1+|x|^2)F_1\|_{H^1}\lesssim \frac{1}{\l^3}\|\tu\|_{L^2}+\frac{1}{\l^2}\|\tu\|_{H^1}\lesssim \sqrt{\alpha^*}\l^{m-2},$$
and thus:
\be\label{b3demi3}
\|(1+|x|^2)F_1\|_{L^2_{[t,T_\eta]}H^1}\lesssim \sqrt{\alpha^*}\l^{m-\frac{3}{2}}.
\ee
The $F_2$ term is estimated from Sobolev embeddings and standard product estimates in Besov spaces:
\bee
\|\nabla^{\frac{3}{2}}F_2\|_{L^{\frac{4}{3}}}&\lesssim& \|\nabla^{\frac{3}{2}}\tu\|_{L^2}(\|\tq\|_{L^4}\|z\|_{L^\infty}+\|z\|_{L^8}^2+\|w\|_{L^8}\|\tu\|_{L^8}+\|\tu\|_{L^8}^2)\\
&&+\|\tu\|_{L^8}(\|\nabla^{\frac{3}{2}}\tq\|_{L^2}\|z\|_{L^8}+\|\tq\|_{L^8}\|\nabla^{\frac{3}{2}}z\|_{L^2}+\|\nabla^{\frac{3}{2}}z\|_{L^2}\|z\|_{L^8}\\
&&+\|\nabla^{\frac{3}{2}}w\|_{L^2}\|\tu\|_{L^8})\\
&\lesssim&\sqrt{\alpha^*} \l^{-\frac{1}{2}}\|\nabla^{\frac{3}{2}}\tu\|_{L^2}+\alpha^*\l^{m-\frac{5}{4}},
\eee
where we used in the last inequality the a priori bound given by \eqref{aprioritu} and \eqref{improvedboundsquareroot} for $\tu$ and the estimate \eqref{cnofeof} for $z$. In turn we obtain:
$$\|\nabla^{\frac{3}{2}}F_2\|_{L^{\frac{4}{3}}_{[t,T_\eta]}L^{\frac{4}{3}}}\lesssim \sqrt{\alpha^*}\l^{\frac{1}{4}}\|\nabla^{\frac{3}{2}}\tu\|_{L^\infty_{[t,T_\eta]}L^2}+\alpha^*\l^{m-\frac{1}{2}}.$$
Injecting this together with \eqref{b3demi2}, \eqref{b3demi3} into \eqref{b3demi1} yields:
$$\|\nabla^{\frac{3}{2}}\tu\|_{L^\infty_{[t,T_\eta]}L^2}\lesssim\sqrt{\alpha^*} \l^{\frac{1}{4}}\|\nabla^{\frac{3}{2}}\tu\|_{L^\infty_{[t,T_\eta]}L^2}+\sqrt{\alpha^*}\l^{m-\frac{3}{2}}$$
and \fref{htrheehalfbound} follows.\\
This concludes the proof \fref{estktun} and of the bootstrap Lemma \ref{bootlemma}.

%%%%%%%%%%%%%%%%%%%%%%%%%%%%%%%%%%%%%%
 %%%%%%%%%%%%%%%%%%%%%%%%%%%%%%%%%%%%%%
 
\section{Instability of Bourgain Wang solutions}
\label{seciofheio}

%%%%%%%%%%%%%%%%%%%%%%%%%%%%%%%%%
%%%%%%%%%%%%%%%%%%%%%%%%%%%%%%%%%%%%%%

We are now in position to prove the main results of the paper which proof we split in several propositions. In the whole section, we let 
\be
\label{hypm}
m\geq 7
\ee
and let 
\be\label{hyptzstar}
z^*\in H^{2m+3}\cap \Sigma
\ee
radially symmetric satisfy \fref{smallnessbis}, \fref{flatattheorigin} with $\alpha^*$, $|\eta|<\eta^*(\alpha^*)$ small enough so that Proposition \ref{propbackwardintegration} holds. We let $z(t)\in \mathcal C(\Bbb R_-,H^{2m+3}\cap \Sigma)$ be the solution to \fref{defutilde}.\\
We claim that the conclusions of Theorem \ref{thmmain} and Theorem \ref{themcontinuation} hold for this class of profiles.

%%%%%%%%%%%%%%%%%%%%%%%%%%%%%%%%%

\subsection{Limit as $\eta\to 0$}

 %%%%%%%%%%%%%%%%%%%%%%%%%%%%%%%%%%%%%%
 %%%%%%%%%%%%%%%%%%%%%%%%%%%%%%%%%%%%%%

We first claim that $u_{\eta}$ constructed by Proposition \ref{propbackwardintegration}  converges strongly as $\eta\to0$ on $(-\infty,0)$ to a Bourgain Wang solution with profile $z^*$.

 \begin{Lemma}[Limit as $\eta\to0$]
 \label{lemmalimit}
Let a sequence 
\be
\label{hyptteta}
T_{\eta}\leq 0, \ \ T_{\eta}\to 0
\ee
such that \fref{dataatteta0} holds. Let $u_{\eta}\in \mathcal C((-\infty,T_{\eta}],H^{\frac 32}\cap \Sigma)$ be the solution to \fref{dataatteta0} constructed in Proposition \ref{propbackwardintegration}. Let $t_1<0$ be the corresponding time\footnote{independent of $\eta$ from Proposition \ref{propbackwardintegration}} of local backwards control. Assume that the phase shift satisfies: 
\be
\label{convergencephaseshidt}
e^{i\gamma^0_{\eta}}\to e^{i\theta}\ \ \mbox{as} \ \ \eta\to 0.
\ee
Then up to a subsequence, $$\forall t<0,\ \ u_\eta(t)\to u_c(t) \ \ \mbox{in} \ \ \Sigma$$ where $u_c\in \mathcal C((-\infty,0),H^{\frac 32}\cap \Sigma)$ is a solution to \fref{nls} which scatters to the left and satisfies the local control: $\forall t\in [t_1,0),$
\be\label{mardi}
\ \ \|u_c(t)-S(t)e^{i\theta}-z(t)\|_{L^2}\lesssim \sqrt{\alpha^*}|t|^{m-2},\,\|u_c(t)-S(t)e^{i\theta}-z(t)\|_{H^1}\lesssim \sqrt{\alpha^*}|t|^{m-3}.
\ee
 \end{Lemma} 
 
  \begin{remark}
 \label{remarkregularite}
For  $\eta>0$ and $T_{\eta}=0$, \fref{dataatteta0} holds from \fref{parametersinitiallyone}. Lemma \ref{lemmalimit} thus yields automatically the existence of a Bourgain solution in the limit $\eta\to 0$ of the solution $u_\eta$ constructed in Proposition \ref{propbackwardintegration} with $\gamma_{\eta}^0=0$, and the obtained solution  $u_c$ is automatically unstable by scattering. A careful track of constants shows that the whole construction requires $m>1$ only and hence a substantial gain on the regularity and more importantly degeneracy in \fref{flatattheorigin}  with respect to the initial Bourgain Wang analysis \cite{BW}. The large $m$ assumption \fref{hypm} is needed first to get uniform bounds on the phase in \fref{uniformbounds}\footnote{where $m>2$ is needed} and more importantly to prove some uniqueness statement about $u_c$, see Proposition \ref{propuniq} below.
 \end{remark}

{\bf Proof of Lemma \ref{lemmalimit}}:\\

{\bf step 1} $\Sigma$ compactness.\\

 We first claim  that $u_\eta(t_1)$ is compact in $\Sigma$ as $\eta\goto 0$. Using Proposition \ref{propbackwardintegration} and in particular \eqref{uniformbounds}, we have the bound:
\be\label{mercre1}
\|u_\eta(t_1)\|_{H^{\frac32}}\lesssim 1,
\ee
which shows that $(u_\eta(t_1))_{0<\eta<\eta^*}$ is compact in $\Sigma(r<R)$ as $\eta\goto 0$ for all $R>0$. The $\Sigma$ compactness of $u_\eta(t_1)$ is now a consequence of a standard localization procedure. Indeed, let a cut off function $\chi(x)=0$ for $|x|\leq 1$ and $\chi(x)=1$ for $|x|\geq 2$, then 
$$\left|\frac{d}{dt}\int \chi_R|u_\eta|^2\right|=2\left|Im\left(\int \nabla\chi_R\cdot\nabla u_\eta)\overline{u_\eta}\right)\right|\lesssim \frac{1}{R},$$ 
$$ \left|\frac{d}{dt}\int \chi_R\left(\frac{1}{2}|\nabla u_\eta|^2-\frac{1}{4}|u_\eta|^4\right)\right|=\left|Im\left(\int \nabla\chi_R\cdot\nabla u_\eta\overline{(\Delta u_\eta+u_\eta|u_\eta|^2)}\right)\right|\lesssim \frac{1}{R}$$ where we used \fref{uniformbounds} and \eqref{cnofeof}. Integrating this backwards from $T_{\eta}$ to $t_1$ and using \fref{dataatteta} yields:
\be\label{montag1}
\lim_{R\goto +\infty}\sup_{0<\eta<\eta^*}\|u_\eta(t_1)\|_{H^1(r>R)}=0. 
\ee
Let now $\psi(x)=0$ for $|x|\leq 1$ and $\psi(x)=|x|^2$ for $|x|\geq 2$, and $\psi_R(x)=R^2\psi(\frac{x}{R})$, $|\psi'_R|^2\lesssim \psi_R$ with constant independent of $R$. Then the decomposition \fref{defqtilde} and the bound \eqref{montag1} ensure:
\bee
 \frac 12\left|\frac{d}{dt}\int \psi_R|u_\eta|^2\right|  & \lesssim & \|w\|_{\Sigma(r>R)}^2+\|\nabla\tu\|_{L^2(r>R)}+\|(\nabla\psi_R)\tu\|_{L^2}^2\\
& \lesssim & o(1)+\int\psi_R|u_{\eta}|^2\textrm { where }o(1)\goto 0\textrm{ as }R\goto 0.
\eee
Integrating this from $t_1$ to $T_{\eta}$ with Gronwall lemma and using \fref{dataatteta} yields:
$$\lim_{R\goto +\infty}\sup_{0<\eta<\eta^*}\| xu_\eta(t_1)\|_{L^2(|x|>R)}=0,$$
which together with \eqref{montag1} and the $\Sigma(r<R)$ compactness of $(u_\eta(t_1))_{0<\eta<\eta^*}$ provided by  \eqref{mercre1} implies up to a subsequence: $$u_\eta(t_1)\to u_c(t_1) \ \ \mbox{in} \ \ \Sigma\ \ \mbox{as} \ \ \eta\to 0.$$ 

{\bf step 2} The limit $u_c$ is a Bourgain Wang solution blowing up at $T=0$.\\

 Let then $u_c\in \mathcal C([t_1,T_c),\Sigma)$ be the solution to \fref{nls} with initial data $u_c(t_1)$, then the $\Sigma$ continuity of the flow ensures: $\forall t\in[t_1,\min(T_c,0))$, $$u_\eta(t)\to  u_c(t) \ \ \mbox{in} \ \ \Sigma.$$ Let $(\lambda(t),b(t),\gamma(t),\e(t))$ be the geometrical decomposition \fref{decompopofupe} associated to $u_{\eta}(t)$ where we have dropped for consistency the $\eta$ dependance, then $u_c$ admits on $[t_1,\min(T_c,0))$ a geometrical decomposition  of the form $$u_c=\frac{1}{\lambda_c(t)}(Q_{0,b_c(t)}+\e_c)\left(t,\frac{x}{\lambda_c(t)}\right)e^{i\gamma_c(t)}e^{i\theta}+z.$$ with: $\forall t\in [t_1,\min(T_c,0))$,
$$\l(t)\goto \l_c(t), \ \ b(t)\goto b_c(t), \ \ \gamma(t)\goto \gamma_c(t), \ \ \e(t)\to \e_c(t)\ \ \mbox{in} \ \ \Sigma$$ as $\eta\to 0$, see \cite{MR3} for related statements. By passing to the limit in \eqref{uniformbounds} and using the explicit formula \fref{parametersinitiallyonebis} and the convergence \fref{convergencephaseshidt}, we obtain the bounds: $\forall t\in [t_1,\min(T_c,0))$,
$$|b_c+t|+|\lambda_c-|t||\lesssim \sqrt{\alpha^*}|t|^m, \ \ \left|\gamma_c+\frac{1}{t}\right|\lesssim\sqrt{\alpha^*} |t|^{m-2},$$
$$\|\e_c\|_{H^1}\lesssim\sqrt{\alpha^*} |t|^{m+1}, \ \ \|\e_c\|_{H^{\frac32}}\lesssim\sqrt{\alpha^*} |t|^{m-\frac32}.$$
This yields that $u_c\in \mathcal C([t_1,0),H^{\frac 32}\cap \Sigma)$ and blows up at $T_c=0$ with $$\|u_c(t)-S(t)e^{i\theta}-z(t)\|_{L^2}\lesssim \sqrt{\alpha^*}|t|^{m-2},\,\|u_c(t)-S(t)e^{i\theta}-z(t)\|_{H^1}\lesssim\sqrt{\alpha^*} |t|^{m-3},$$ and \fref{mardi} is proved. Now $u_c$ is global to the left and scatters for $\alpha^*$ small enough from \fref{mardi} and Lemma \ref{lemmaopen}. This concludes the proof of Lemma \ref{lemmalimit}.

 %%%%%%%%%%%%%%%%%%%%%%%%%%%%%%%%%%%%%%
 %%%%%%%%%%%%%%%%%%%%%%%%%%%%%%%%%%%%%%

 \subsection{Conditional uniqueness of the Bourgain Wang solutions}

 %%%%%%%%%%%%%%%%%%%%%%%%%%%%%%%%%%%%%%
 %%%%%%%%%%%%%%%%%%%%%%%%%%%%%%%%%%%%%%
 
 We now claim the following conditional uniqueness statement about the Bourgain Wang type solutions satisfying \fref{mardi}:
 
 \begin{Prop}[Conditional uniqueness of the Bourgain Wang solutions]
 \label{propuniq}
Let $z^*\in H^{2m+3}$ a radially symmetric function satisfying \eqref{flatattheorigin} with $m=7$. 
Let $t_1<0$. Then, there is a unique $u\in \mathcal C([t_1,0),H^{\frac 32})$ solution to \fref{nls} with 
\be\label{uniq}
\|u(t)-S(t)-z(t)\|_{L^2}\lesssim |t|^5,\,\|u(t)-S(t)-z(t)\|_{H^1}\lesssim |t|^4.
\ee
 \end{Prop}

This proposition follows by a simple further use of the Lyapounov control of Lemma \ref{lemma:timederivative} as in \cite{RS}. The proof is postponed to Appendix A.

%%%%%%%%%%%%%%%%%%%%%%%%%%%%%%%%%%%%%%
 %%%%%%%%%%%%%%%%%%%%%%%%%%%%%%%%%%%%%%
 
 \subsection{Proof of the main theorems}
 
 %%%%%%%%%%%%%%%%%%%%%%%%%%%%%%%%%%%%%%
 %%%%%%%%%%%%%%%%%%%%%%%%%%%%%%%%%%%%%%
 
 We are now in position to conclude the proof of the main Theorems.\\
 
 {\bf Proof of Theorem \ref{thmmain}, Theorem \ref{themcontinuation}}\\
 
{\bf step 1} Definition and continuity of the map $\Gamma$.\\

 Let $m\geq 7$ and $z^*$ radially symmetric satisfy \fref{hyptzstar}, \fref{smallnessbis}, \fref{flatattheorigin} with $\alpha^*$, $|\eta|<\eta^*(\alpha^*)$  small enough. Let $0<\delta\ll1$ a small enough universal number\footnote{depending only on $z^*$, see \fref{poutfmoins}} to be chosen later. For $\eta\in (-\eta^*,\eta^*)$, $\eta \neq 0$ we let 
 \be
 \label{normalizationtimephase}
 T_{\eta}=\left\{\begin{array}{ll} 0 \ \ \mbox{for} \ \ \eta>0\\ -\sqrt{(1+\delta)|\eta|} \ \ \mbox{for} \ \ \eta<0\end{array}\right ., \ \ \gamma_\eta^0=0.
 \ee 
 Observe that \fref{dataatteta0} is fulfilled both for $\eta>0$ and $\eta<0$ from the explicit law \fref{parametersinitiallyone}. We then let $u_{\eta}\in \mathcal C((-\infty,T_{\eta}],H^{\frac32}\cap \Sigma)$ be the solution to \fref{nls} constructed in Proposition \ref{propbackwardintegration} and define the map:
 \be
 \label{defgamma}
 \Gamma(\eta)=u_{\eta}(-1), \ \ \eta\neq 0
 \ee
For $\eta=0$, the uniqueness statement of Proposition \ref{propuniq} together with the compactness statement of Lemma \ref{lemmalimit} ensures: $$\forall t<0, \ \ u_{\eta}(t)\to u_{BW}^0(t) \ \ \mbox{in} \ \ H^1 \ \ \mbox{as} \ \ \eta\to 0$$ where $u_{BW}^0$ is the unique Bourgain Wang solution with regular profile $z^*$ and singular profile $S(t)$ as $t\to 0$ satisfying the bounds \fref{uniq}. We thus define $$\Gamma(0)=u_{BW}^0(-1)$$ so that $$\Gamma:(-\eta^*,\eta^*)\to  H^1\ \ \mbox{is continuous}$$ and $\Gamma(\eta)\in H^{\frac32}\cap \Sigma$. It remains to show the expected behavior of $u_{\eta}(t)$ for $\eta\neq 0$.\\

{\bf step 2} Scattering and continuation after blow up time.\\

Let $\eta>0$, then $u_{\eta}\in \mathcal C((-\infty,0],H^{\frac 32}\cap \Sigma)$ from Proposition \ref{propbackwardintegration}. Let $v_{\eta}$ be the solution to \fref{nls} with initial data 
\be
\label{initialdata}
v_{\eta}(0)=\overline{u_\eta(0)}
\ee
then from \fref{solutionexplicite}, \fref{parametersinitiallyone} and Lemma \ref{scatteringprofiles}, we have:
$$v_{\eta}(0)=\frac{1}{\sqrt{\eta}}P_\eta\left(\frac{x}{\sqrt{\eta}}\right)e^{i\gamma_\eta(0)}+\overline{z^*}=\tilde{Q}_\eta(0)e^{i\gamma_{\eta}^0}+\overline{z^*}
$$
with 
\be
\label{defphaseshift}
\gamma^0_{\eta}=-2\int_{-1}^{0}\frac{d\tau}{\eta+\tau^2}.
\ee
The profile $\overline{z^*}$ satisfies the assumptions of Proposition \ref{propbackwardintegration}, and thus $v_{\eta}\in \mathcal C((-\infty,0],H^{\frac 32}\cap \Sigma)$ and scatters to the left. Now the time reversibility of the (NLS) flow ensures 
\be
\label{forulaionevrsebarre}
\forall t>0, \ \ u_{\eta}(t)=\overline{v_{\eta}(-t)}
\ee and thus $u_{\eta}$ is global and scatters to the right. Moreover, let $\theta\in \Bbb R$. From \fref{defphaseshift}, we have:
$$\gamma^0_\eta \goto-\infty\textrm{ when }\eta\goto 0_+.$$
Thus, there exists a sequence $\eta_n\to 0$ such that $$e^{i\gamma_{\eta_n}^0}\to e^{i\theta} \ \ \mbox{as} \ \ n\to \infty.$$ Along this sequence, we apply Lemma \ref{lemmalimit} to $v_{\eta_n}$ which together with the uniqueness statement of Proposition \ref{propuniq} ensures:
$$\forall \ \ t<0,\  \ v_{\eta_n}(t)\to u_{BW}^{\theta} \textrm{ in }\Sigma\ \ \mbox{as} \ \ \eta\to 0_+,$$  and \fref{opjepoueu} now follows from \fref{forulaionevrsebarre}. This concludes the proof of the case $\eta>0$ in Theorem \ref{thmmain}  and of Theorem \ref{themcontinuation}.\\

{\bf step 3} Sufficient criterion for loglog blow up.\\

We now turn to the case $\eta<0$. Let us start with recalling the following criterion of log-log blow up which follows from \cite{R1}. We let $$Q_{b}=Q_{\eta=0,b}=Qe^{-i\frac{b|y|^2}{4}}.$$
 
 \begin{Prop}[Open characterization of the log-log set, \cite{R1}]
 \label{proploglog}
 Let $\alpha^*>0$ denote a small enough universal constant. Let an initial data of the form $$u_0=\frac{1}{\lambda_0}(Q_{b_0}+\e_0)\left(\frac{x}{\lambda_0}\right)e^{i\gamma_0}$$ where $\e_0\in H^1_{rad}$ satisfies the orthogonality conditions:
 \be
 \label{orthoieoglogprop}
 Re(\e_0,\overline{|y|^2 Q})=Im(\e_0,\overline{\Lambda Q})=Im(\e_0,\overline{\Lambda^2Q})
 \ee
 and assume that the following bounds hold:\\
 (i) $L^2$ control: 
 \be
 \label{ltwocontorl}
 \int Q^2<\int |u_0|^2\leq \int Q^2+\alpha^*;
 \ee
 (ii) Energy control: 
 \be
 \label{ebergyintialdtat}
 E(u)\leq \alpha^*\int|\nabla u_0|^2;
 \ee
 (iii) Open characterization of the log log set: 
\be\label{mercredi3}
 f_-(0)=\frac{b_0}{\lambda_0}-\frac{\sqrt 8}{\|yQ\|_{L^2}}\sqrt{|E_0|}>0,
 \ee
 then $u$ blows up in finite time $T>0$ in the log log regime with the a priori bound: 
 \be
 \label{upperboundt}
 T\lesssim \frac{\lambda_0}{f_-(0)}.
 \ee
 \end{Prop}
 
 {\bf Proof of Proposition \ref{proploglog}}: We apply Lemma 6 in \cite{R1} which ensures that $u$ blows up in finite time $T>0$ in the log log regime. It remains to prove the upper bound \fref{upperboundt}. Following the notations of the proof of Lemma 6 in \cite{R1}, $u$ admits on $[0,T)$ a decomposition $$u(t,x)=\frac{1}{\lambda}(Q_{b(t)}+\e)\left(t,\frac{x}{\lambda(t)}\right)e^{i\gamma(t)}$$ where $\e$ satisfies the orthogonality conditions \fref{orthoieoglogprop}. The scaling parameter satisfies the following estimate:
$$\forall t\in [0,T),   \ \  \lambda(t)\leq 2\lambda_0,$$ and there holds the rigidity:$$\forall t\in [0,T), \ \ f_-(t)=\frac{b(t)}{\lambda(t)}-\frac{\sqrt 8}{\|yQ\|_{L^2}}\sqrt{|E_0|}>0.$$ More precisely, $f_-$ satisfies the differential inequality: $$-\frac{d}{dt}\left(\frac{1}{f_-}\right)\geq \frac{C}{\l}\gtrsim \frac{1}{\lambda_0},$$ see p599 of \cite{R1}. Integrating this from $0$ to $T$ yields the bound: $$\frac{T}{\lambda_0}\lesssim \frac{1}{f_-(0)},$$ this is \fref{upperboundt}. This concludes the proof of Proposition \ref{proploglog}.\\

{\bf step 4} Loglog blow up for $\eta<0$.\\

Let now $\eta<0$ and $u_{\eta}\equiv u$ be the solution to \fref{nls} given by Proposition \ref{propbackwardintegration} with the normalization \fref{normalizationtimephase}. Then $u\in \mathcal C((-\infty,T_{\eta}],H^{\frac32}\cap \Sigma)$ and scatters to the left. Let us check that $u$ satisfies the assumptions \fref{ltwocontorl}, \fref{ebergyintialdtat}, \fref{mercredi3}
 of Proposition \ref{proploglog} which will yield the claim.\\
 
We compute the energy of $u$. Recall from \eqref{parametersinitiallyone} that 
\be
\label{vnovhrohv}
\lambda_\eta(T_{\eta})=\sqrt{\delta |\eta|}, \ \ b_{\eta}(T_{\eta})=\sqrt{(1+\delta)|\eta|}.
\ee
We then compute from \eqref{uniformdeacy}, \fref{coputationngry}, \fref{defqetab}:
$$E(\qte(T_{\eta}))=\frac{\|yQ\|_{L^2}^2}{8\l_\eta^2(T_{\eta})}\left[b^2_{\eta}(T_{\eta})-|\eta|+o(\eta)\right]=\frac{\|yQ\|_{L^2}^2}{8}+o(1) \ \ \mbox{as} \ \ \eta\to 0.$$
We thus estimate from \fref{bins2}: 
\bea
\nonumber
E(u)&=&E(\qte(T_{\eta}))+E(z^*)+O\left(\frac{1}{\l^2}(\|\e_ze^{-c|y|}\|_{L^\infty}+\|\e_ze^{-c|y|}\|_{L^\infty}^3)\right)\\
\nonumber
&=&E(\qte(T_{\eta}))+E(z^*)+O(\l_\eta^{m}(T_{\eta}))\\
\label{cnekoeoeu}&=&\frac{\|yQ\|_{L^2}^2}{8}+E(z^*)+o(1) \ \ \mbox{as} \ \ \eta\to 0.
\eea
Let us now introduce from standard argument the unique decomposition 
\be\label{mercre3}
u(T_{\eta})=\frac{1}{\lambda}(Q_{b}+\e)\left(\frac{x}{\lambda}\right)e^{i\gamma}=\frac{1}{\lambda_{\eta}(T_{\eta})}(Q_{\eta,b_{\eta}(T_{\eta})}+\e_z(T_{\eta}))\left(\frac{x}{\lambda_{\eta}(T_{\eta})}\right)e^{i\gamma_{\eta}(T_{\eta})}
\ee
with $\e$ satisfying the orthogonality conditions \fref{orthoieoglogprop}. Then, taking the scalar product of \eqref{mercre3} with the three orthogonality conditions \eqref{orthoieoglogprop} yields:
$$\left|\frac{\lambda_\eta(T_\eta)}{\lambda}-1\right|+|b_{\eta}(T_\eta)-b|+|\gamma-\gamma_{\eta}(T_\eta)|\lesssim  \|Q_b-Q_{\eta,b_\eta(T_\eta)}\|_{L^2}+\|e^{-c|y|}\e_z(T_\eta)\|_{L^2},$$
which together with \eqref{uniformdeacy} and \eqref{bins2} ensures the bounds: 
\be
\label{boundss}
|b-b_{\eta}(T_{\eta})|+\left|\frac{\lambda}{\lambda_{\eta}(T_{\eta})}-1\right|+|\gamma-\gamma_{\eta}(T_{\eta})|\lesssim |\eta|.
\ee
 Hence:
\bee
\frac{b}{\l} & = & \frac{b_{\eta}(T_{\eta})}{\lambda_{\eta}(T_{\eta})}+O\left(\frac{|b-b_{\eta}(T_{\eta})|}{\lambda_{\eta}(T_{\eta})}+\frac{|b_{\eta}(T_{\eta})|}{\lambda_{\eta}(T_\eta)}\left|\frac{\lambda}{\lambda_\eta(T_{\eta})}-1\right|\right)\\
& = & \sqrt{1+\frac{1}{\delta}}(1+o(1)) \ \mbox{as} \ \ \eta\to 0.
\eee
Together with \fref{cnekoeoeu}, this ensures that for a given $z^*$, we may chose $\delta>0$ small enough and find $\eta^*>0$ such that for all $-\eta^*<\eta<0$, 
\be
\label{poutfmoins}
f_-(T_{\eta})=\frac{b}{\lambda}-\frac{\sqrt 8}{\|yQ\|_{L^2}}\sqrt{E_0}\geq \frac{1}{2\sqrt{\delta}}>0,
\ee
which is \eqref{mercredi3}. 

$\delta$ being now fixed, we have from \fref{vnovhrohv}, \fref{cnekoeoeu} and the fact that $\|\nabla u(T_\eta)\|_{L^2}\sim \l_\eta(T_\eta)^{-1}$:
$$\frac{E(u)}{\|\nabla u(T_\eta)\|^2_{L^2}}\lesssim \l^2_\eta(T_\eta)\lesssim |\eta|$$
and \fref{ebergyintialdtat} follows  for $|\eta|<\eta^*(\a^*)$ small enough. As for the $L^2$ control \fref{ltwocontorl}:
\bee
\int |u|^2&=&\int |\qte(T_\eta)|^2+2\Re\left(\int z(T_\eta)\overline{\qte(T_\eta)}\right)+\int |z(T_\eta)|^2\\
&=& \int |P_\eta|^2+\int |z^*|^2+O\left(\|\e_ze^{-c|y|}\|_{L^2}\right)\\
&=& \int Q^2+\frac{|\eta|}{4}\|yQ\|^2_{L^2}+\int |z^*|^2+O(\eta^2),
\eee
where we used \eqref{invariantl2peta}, \eqref{bins2} and \eqref{vnovhrohv} in the last inequality. Hence 
\fref{ltwocontorl} holds for $|\eta|<\eta^*(\a^*)$ small enough.

Finally, the assumptions \eqref{ltwocontorl}-\eqref{mercredi3} hold. Hence , we conclude from Proposition \ref{proploglog} that $u$ blows up in the log log regime in forward time at some time $T_{\eta}^*>T_{\eta}$ with from \fref{upperboundt}, \fref{boundss}, \fref{poutfmoins}:$$T_{\eta}^*\leq T_{\eta}+\frac{\lambda}{f_-(T_{\eta})}\leq-\sqrt{(1+\delta)|\eta|}+C\delta\sqrt{|\eta|}+O(\eta)\leq -\frac{\sqrt{|\eta|}}{2}.$$ 
This concludes the proof of Theorem \ref{thmmain}.\\

{\bf Proof of Corollary \ref{corforward}}: Let $u_\eta=\Gamma(\eta)$ be the solution to \fref{nls} given by Theorem \ref{thmmain}. For $\eta<0$, let $T_{\eta}^*<0$ be the blow up time of $u_{\eta}$, and $T_{\eta}^*=0$ for $\eta\geq 0$. Note from \cite{MR3} that the blow up time in the loglog regime  is a continuous function of the initial data in $H^1$ and hence the map $\eta\mapsto T_{\eta}^*$ is continuous. We then define the transformations: 
$$v^1_{\eta}(1,x)=\Gamma^1(\eta)=u_{\eta}(-1,x)e^{i\frac{|x|^2}{4}}, \ \ v^2_{\eta}(1,x)=\Gamma^2(\eta)=u_{\eta}(T^*_{\eta}-1,x)e^{i\frac{|x|^2}{4}}$$ which are continuous maps from $[-1,1]\to \Sigma$. From the pseudo conformal invariance, we have the explicit deformation: $$\forall \tau>0,  \ \ v^1_{\eta}(\tau,x)=\frac{1}{\tau}u_{\eta}\left(\frac{-1}{\tau},\frac{x}{\tau}\right)e^{i\frac{|x|^2}{4\tau}},\ \ v^2_{\eta}(\tau,x)=\frac{1}{\tau}u_{\eta}\left(T^*_{\eta}-\frac{1}{\tau},\frac{x}{\tau}\right)e^{i\frac{|x|^2}{4\tau}}.$$ The conclusions of Corollary  \ref{corforward} now follow from a direct inspection which is left to the reader. This concludes the proof of Corollary \ref{corforward}.

%%%%%%%%%%%%%%%%%%%%%%%%%%%%%%%%%
%%%%%%%%%%%%%%%%%%%%%%%%%%%%%%%%%

\section*{Appendix A}

%%%%%%%%%%%%%%%%%%%%%%%%%%%%%%%%%
%%%%%%%%%%%%%%%%%%%%%%%%%%%%%%%%%

This appendix is devoted to the proof of the Proposition \ref{propuniq}.\\
Let $u_c$ be the Bourgain Wang type solution constructed in Lemma \ref{lemmalimit} with $\theta=0$. Note that $u_c$ satisfies \fref{uniq} from \fref{mardi} with $m=7$. Let $u$ be another solution satisfying \fref{uniq}. We need to show that $u\equiv u_c$.

%%%%%%%%%%%%%%%%%%%%%%%%%%%%%%%%%%%%%%%%%%%%%%%
 %%%%%%%%%%%%%%%%%%%%%%%%%%%%%%%%%%%%%%%%%%%%%%%

 \subsection{Energy estimates for the flow near $u_c$}
 
 %%%%%%%%%%%%%%%%%%%%%%%%%%%%%%%%%%%%%%%%%%%%%%%
 %%%%%%%%%%%%%%%%%%%%%%%%%%%%%%%%%%%%%%%%%%%%%%%

First, recall from the proof of Lemma \ref{lemmalimit} that $u_c$ admits a geometrical decomposition  of the form $$u_c=\tq_c+z+\tu_c=\frac{1}{\lambda_c(t)}(Q_{b_c(t)}+\e_z+\e_c)\left(s,\frac{x}{\lambda_c(t)}\right)e^{i\gamma_c(t)}, \ \ \frac{ds}{dt}=\frac{1}{\lambda_c^2}$$ with\footnote{the inequalities \eqref{boundpmc}-\eqref{bounduc} follow from the proof of Lemma \ref{lemmalimit} in the case $m=4$ for the vanishing of $z^*$ at the origin. Since this is enough for the uniqueness part of the proof of Proposition \ref{propuniq}, we have chosen to state these inequality with $m=4$ instead of 
$m=7$}:
\be\label{boundpmc}
|b_c+t|+|\lambda_c-|t||\lesssim |t|^4, \ \ \left|\gamma_c+\frac{1}{t}\right|\lesssim |t|^2, 
\ee
\be\label{boundmodc}
|Mod_c(t)|=|(b_c)_s+b_c^2|+\left|\frac{(\l_c)_s}{\l_c}+b_c\right|+|(\tilde{\gamma}_c)_s|\lesssim |t|^5, 
\ee
\be\label{boundec}
\|\e_c\|_{H^1}\lesssim |t|^5, \ \ \|\e_c\|_{H^{\frac32}}\lesssim |t|^{\frac52},
\ee
and
\be\label{bounduc}
 \|\tu_c\|_{L^2}\lesssim |t|^5, \ \  \|\tu_c\|_{H^1}\lesssim |t|^4, \ \ \|\tu_c\|_{H^{\frac32}}\lesssim |t|.
\ee

Let us now decompose:
\be
\label{un50}
u=u_c+\ttu,  \ \ \ds\ttu(t,x)=\frac{1}{\l_c(t)}\e\left(t,\frac{x}{\l_c(t)}\right)e^{i\gamma_c(t)}.
\ee
Here we do not impose modulation theory and orthogonality conditions on $\e$. We however claim that the a priori estimate from \eqref{uniq}:
\be
\label{backwardttu}
 \|\ttu\|_{L^2}\lesssim |t|^5,\, \|\ttu\|_{H^1}\lesssim |t|^4
 \ee
 is enough to treat the instability generated by the null space of $L$ perturbatively.\\
 Let 
 \be\label{un70}
N(t):=\sup_{t<\tau<0}\left(\|\ttu(\tau)\|^2_{H^1}+\frac{\|\ttu(\tau)\|^2_{L^2}}{\l_c(\tau)^2}\right),
\ee 
 and 
 \be
 \label{defst}
 \mbox{Scal}(t)  =  (\e_1,Q)^2+(\e_2,\Lambda Q)^2+(\e_1,|y|^2Q)^2+(\e_2,\rho)^2.
 \ee
 
 We first claim the following energy bound:
 
 \begin{Lemma}
 \label{lemma:basicestimateun}
 There holds for $t$ close enough to 0: 
 \be
 \label{controlnscal}
 N(t)\lesssim \sup_{t\leq \tau<0}\frac{\mbox{Scal}(\tau)}{\l_c(\tau)^2}+\int_t^0\frac{\mbox{Scal}(\tau)}{\l_c(\tau)^3}d\tau.
 \ee
 \end{Lemma}
 
 {\bf Proof of Lemma \ref{lemma:basicestimateun}}\\
 
 It is a consequence of the energy estimate \fref{crc6} together with the a priori bound \fref{bounduc}.\\
 
 {\bf step 1} Application of Lemma \ref{lemma:timederivative}.\\
 
 Let $$w=u_c=(u_c)_1+i(u_c)_2$$ and $\mathcal I(\ttu)$ be given by \fref{defI}, we claim that:
 \be
 \label{iegigi}
\frac{\|\ttu\|_{L^2}^2}{\lambda_c^3}+O\left(N(t)+\frac{\mbox{Scal}(t)}{\lambda_c^3}\right)\lesssim  \frac{d\mathcal I}{dt}.
\ee 
 Indeed, we apply Lemma \ref{lemma:timederivative} with $w=u_c$, then the bound \fref{aprioriboundwsw} holds from \eqref{cnofeof} \fref{bounduc}, and $\psi$ given by \fref{eqwgobale} 
 is identically zero. Furthermore, the bounds \eqref{aprioritubis} \eqref{roughboundmodulationbis} hold from \eqref{boundmodc} \eqref{bounduc}. Hence \fref{crc6} becomes:
 \bea
 \label{un55}
\nonumber & & \frac{d\mathcal I}{dt}= -\frac{1}{\l_c^2}\Im\left(\int u_c^2\overline{\ttu}^2\right)-\Re\left(\int \p_tu_c\overline{(2|\ttu|^2u_c+\ttu^2\overline{u_c})}\right)\\
 \nonumber & + & \frac{b_c}{\l_c^2}\Bigg\{\int \frac{|\ttu|^2}{\l_c^2}+\Re\left(\int\nabla^2\phi\left(\frac{x}{A\l_c}\right)(\n\ttu,\overline{\n\ttu})\right) -\frac{1}{4A^2}\left(\int\Delta^2\phi\left(\frac{x}{\l_c}\right)\frac{|\ttu|^2}{\l_c^2}\right)\\
\nonumber & + & \l_c\Re\left(\int A\nabla\phi\left(\frac{x}{A\l_c}\right)(2|\ttu|^2u_c+\ttu^2\overline{u_c})\cdot\overline{\nabla u_c}\right)\Bigg\}\\
& + & O\left(\frac{\|\ttu\|^{2}_{L^2}}{\l_c^2}+\|\ttu\|^{2}_{H^1}\right).
\eea
We consider the first two terms in the right-hand side of \eqref{un55} and expand $u_c=w_c+\tu_c$, where $w_c=\tq_c+z$:\footnote{Keep in mind that we do not have satisfactory well localized bounds in $\tilde{u}_c$, and the corresponding terms will be treated using the smallness \fref{bounduc}}
\bea
\label{un56}
& & -\frac{1}{\l_c^2}\Im\left(\int u_c^2\overline{\ttu}^2\right)-\Re\left(\int \p_tu_c\overline{(2|\ttu|^2u_c+\ttu^2\overline{u_c})}\right)\\
\nonumber & = & -\frac{1}{\l_c^2}\Im\left(\int w_c^2\overline{\ttu}^2\right)-\Re\left(\int \p_tw_c\overline{(2|\ttu|^2w_c+\ttu^2\overline{w_c})}\right)\\
\nonumber & & -\frac{1}{\l_c^2}\Im\left(\int (2\tu_cw_c+\tu_c^2)\overline{\ttu}^2\right)-\Re\left(\int \p_tw_c\overline{(2|\ttu|^2\tu_c+\ttu^2\overline{\tu_c})}\right)\\
\nonumber & - & \Re\left(\int \p_t\tu_c\overline{(2|\ttu|^2w_c+\ttu^2\overline{w_c})}\right)-\Re\left(\int \p_t\tu_c\overline{(2|\ttu|^2\tu_c+\ttu^2\overline{\tu_c})}\right).
\eea
Arguing like for the proof of \eqref{firtboundq}, we may rewrite the first two terms in the right-hand side of \eqref{un56} as:
\bea
\label{un57}
& & -\frac{1}{\l_c^2}\Im\left(\int w_c^2\overline{\ttu}^2\right)-\Re\left(\int \p_tw_c\overline{(2|\ttu|^2w_c+\ttu^2\overline{w_c})}\right)\\
\nonumber & = &  -\frac{b_c}{\l^2_c}\int ((|\tq_c|^2+2\ts_c^2)\ttu_1^2+4\ts_c\tt_c\ttu_1\ttu_2+(|\tq_c|^2+2\tt_c^2)\ttu_2^2)\\
\nonumber &- & \frac{b_c}{\l_c}\Re\left(\int \left(\frac{x}{\l_c}\right) (2|\ttu|^2\tq_c+\ttu^2\overline{\tq_c})\cdot\overline{\nabla \tq_c}\right)+  O\left(\frac{\|\ttu\|^{2}_{L^2}}{\l^2_c}+\|\ttu\|^{2}_{H^1}\right)
\eea
where $\tq_c=\ts_c+i\tt_c$. For the next two terms in the right-hand side of \eqref{un56}, we use Sobolev embeddings and \fref{bounduc}  to obtain:
\bea
\label{un58}
\nonumber & & \bigg|-\frac{1}{\l_c^2}\Im\left(\int (2\tu_cw_c+\tu_c^2)\overline{\ttu}^2\right)-\Re\left(\int \p_tw_c\overline{(2|\ttu|^2\tu_c+\ttu^2\overline{\tu_c})}\right)\bigg|\\
\nonumber & \lesssim & \frac{1}{\l_c^2}\|w_c\|_{L^\infty}\|\tu_c\|_{L^2}\|\ttu\|^2_{L^4}
+\frac{1}{\l_c^2}\|\tu_c\|^2_{L^4}\|\ttu\|^2_{L^4}\\
& + & \|\p_tw_c\|_{L^\infty}\|\tu_c\|_{L^2}\|\ttu\|^2_{L^4}\lesssim  \frac{\|\ttu\|^{2}_{L^2}}{\l^2_c}+\|\ttu\|^{2}_{H^1}
\eea
where we used the bound\footnote{The worst term is generated by the phase $|(\gamma_c)_t|\lesssim \frac{1}{\lambda_c^2}$.} $$\|\pa_tw_c\|_{L^{\infty}}\lesssim \frac{1}{\lambda_c^3}.$$
The last two terms in the right-hand side of \eqref{un56} require using the equation satisfied by $\tu_c$:
\be\label{tintin}
i\p_t\tu_c=-\Delta\tu_c-(|u_c|^2u_c-|w_c|^2w_c)-\psi_c
\ee
where $\psi_c$ is defined by:
$$\psi_c=i\p_t w_c+\Delta w_c+|w_c|^2w_c.
$$
Recall from \eqref{estpsi}:
$$\|\psi_c\|_{L^2}\lesssim \frac{|Mod_c(t)|+\l_c^6}{\lambda_c^2}\lesssim \lambda_c^3.$$ Using this together with \eqref{tintin}, integration by parts, Sobolev embeddings and the $H^{\frac{3}{2}}$ bound \fref{bounduc} now yields:
\bea
\label{un58bis}
\nonumber & & \bigg|-\Re\left(\int \p_t\tu_c\overline{(2|\ttu|^2w_c+\ttu^2\overline{w_c})}\right)-\Re\left(\int \p_t\tu_c\overline{(2|\ttu|^2\tu_c+\ttu^2\overline{\tu_c})}\right)\bigg|\\
\nonumber & \lesssim & \|\tu_c\|_{H^{\frac{3}{2}}}\left[\|2|\ttu|^2w_c+\ttu^2\overline{w_c}\|_{H^{\frac{1}{2}}}+\|2|\ttu|^2\tu_c+\ttu^2\overline{\tu_c}\|_{H^{\frac{1}{2}}}\right]\\
\nonumber & + & \|(|u_c|^2u_c-|w_c|^2w_c)+\psi_c\|_{L^2}\left[\|w_c\|_{L^\infty}\|\ttu\|^2_{L^4}+\|\tu_c\|_{L^6}\|\ttu\|^2_{L^6}\right]\\
&\lesssim & \frac{\|\ttu\|^{2}_{L^2}}{\l^2_c}+\|\ttu\|^{2}_{H^1}.
\eea
We now consider the last term in the right-hand side of \eqref{un55} and compute:
\bea
\label{un60}
& &  \Re\left(\int A\nabla\phi\left(\frac{x}{A\l_c}\right)(2|\ttu|^2u_c+\ttu^2\overline{u_c})\overline{\nabla u_c}\right)\\
 \nonumber &= & \Re\left(\int A\nabla\phi\left(\frac{x}{A\l_c}\right)(2|\ttu|^2\tq_c+\ttu^2\overline{\tq_c})\overline{\nabla \tq_c}\right)+\mbox{Error}
 \eea
 with from Sobolev embeddings, \eqref{cnofeof}, \eqref{bins2}, and \fref{bounduc}:
\bea
\label{un61}
\nonumber &  & |\mbox{Error}|\lesssim  \left|\Re\int A\nabla\phi\left(\frac{x}{A\l_c}\right)(2|\ttu|^2(z+\tu_c)+\ttu^2\overline{(z+\tu_c)})\overline{\nabla \tq_c}\right|\\
\nonumber & +&\left|\Re\int A\nabla\phi\left(\frac{x}{A\l_c}\right)(2|\ttu|^2\tq_c+\ttu^2\overline{\tq_c})\overline{\nabla (z+\tu_c)}\right|\\
\nonumber & +& \left|\Re\int A\nabla\phi\left(\frac{x}{A\l_c}\right)(2|\ttu|^2(z+\tu_c)+\ttu^2\overline{(z+\tu_c)})\overline{\nabla(z+\tu_c)}\right|\\
\nonumber & \lesssim & \left(\frac{1}{\l_c^2}\|e^{-c|y|}\e_z\|_{L^2}+\|\n\tq_c\|_{L^\infty}\|\tu_c\|_{L^2}+\|\tq_c\|_{L^\infty}(\|\n z\|_{L^2}+\|\n\tu_c\|_{L^2})\right)\|\ttu\|^2_{L^4}\\
&+ & (\|z\|_{L^6}+\|\tu_c\|_{L^6})(\|\n z\|_{L^2}+\|\n\tu_c\|_{L^2})\|\ttu\|^2_{L^6}
%&\lesssim & (\frac{1}{\l_c^2}\l_c^{4}+\frac{1}{\l_c}\l_c^{3})\|\ttu\|_{L^2(\RR^2)}\|\ttu\|_{H^1(\RR^2)}\\
%&+ &\l_c^{10/3}\l_c^{3}\|\ttu\|^{2/3}_{L^2(\RR^2)}\|\ttu\|^{4/3}_{H^1(\RR^2)}\\
 \lesssim  \frac{\|\ttu\|^{2}_{L^2}}{\l^2_c}+\|\ttu\|^{2}_{H^1}.
\eea
Collecting the estimated \eqref{un55}-\eqref{un61} yields:
\bee
\nonumber & & \frac{d\mathcal I}{dt}  =  \frac{b_c}{\l_c^2}\left[\int \frac{|\ttu|^2}{\l_c^2}+\Re\left(\int\nabla^2\phi\left(\frac{x}{A\l_c}\right)(\n\ttu,\overline{\n\ttu})\right)\right .\\
\nonumber& -& \left . \int ((|\tq_c|^2+2\ts_c^2)\ttu_1^2+4\ts_c\tt_c\ttu_1\ttu_2+(|\tq_c|^2+2\tt_c^2)\ttu_2^2)-\frac{1}{4A^2}\int\Delta^2\phi\left(\frac{x}{A\l_c}\right)\frac{|\ttu|^2}{\l_c^2}\right .\\
\nonumber&+ &\left . \l_c\Re\left(\int \left(A\nabla\phi\left(\frac{x}{A\l_c}\right)-\left(\frac{x}{\l_c}\right)\right)k(x)(2|\ttu|^2\tq_c+\ttu^2\overline{\tq_c})\cdot\overline{\nabla \tq_c}\right)\right]\\
& + & O\left(\frac{\|\ttu\|^{2}_{L^2}}{\l_c^2}+\|\ttu\|^{2}_{H^1}\right).
\eee
We now use the uniform proximity of $Q_{b_c}$ to $Q$ and the coercitivity property \fref{coerclinearenergy} to conclude like for the proof of \fref{firtboundq}:
\bee
 \frac{b_c}{\lambda_c^4}\left[\int|\nabla \e|^2e^{-\frac{|y|}{\sqrt{A}}}+\int|\e|^2+O(\mbox{Scal}(t))\right]+O
\left(\frac{\|\ttu\|^{2}_{L^2}}{\l_c^2}+\|\ttu\|^{2}_{H^1}\right)\lesssim \frac{d\mathcal I}{dt}
\eee
which implies \fref{iegigi}.\\

{\bf step 2} Coercivity of $\mathcal I$.\\

We now recall from \fref{defI} the formula:
\bee
\mathcal I(t) & = & \frac{1}{2}\int |\n\ttu|^2 +\half\int \frac{|\ttu|^2}{\l_c^2}-\frac{1}{4}\int |u_c+\ttu|^4+\frac{1}{4}\int |u_c|^4\\
\nonumber & + & \int |u_c|^2(u_c)_1\ttu_1+\int |u_c|^2(u_c)_2\ttu_2+\half\frac{b_c}{\l_c}\Im\left(\int A\nabla\phi\left(\frac{x}{A\l_c}\right)\cdot\nabla\ttu\overline{\ttu}\right).
\eee
Expanding $u_c=w_c+\tu_c$ and arguing like for the proof of \fref{iegigi}, we get using \fref{backwardttu} the rough upper bound: 
\be
\label{boundarytem}
|\mathcal I|\lesssim \|\ttu(t)\|^2_{H^1}+\frac{\|\ttu(t)\|^2_{L^2}}{\l^2_c(t)}\to 0 \ \ \mbox{as} \ \ t\to 0.
\ee
 We now claim the lower bound:
\bea
\label{coervifityagain}
\nonumber \nonumber \mathcal I(t) & \geq & \frac{1}{2\lambda_c^2} \left[(L_+\e_1,\e_1)+(L_-\e_2,\e_2)+o\left(\|\e\|_{H^1}^2\right)\right] \\
& \geq & \frac{\underline{c}}{\lambda_c^2}\left[\int\|\e\|_{H^1}^2-\mbox{Scal}(t)\right].
\eea
The proof is very similar to the one of \fref{esticoineo} using also the control of interaction terms similar to \fref{un58}, \fref{un58bis}, \fref{un61}. The details are left to the reader.\\
Integrating \fref{iegigi} from $t$ to $0$ using the boundary condition \fref{boundarytem} and the lower bound \fref{coervifityagain} now yields \fref{controlnscal}.\\
This concludes the proof of Lemma \ref{lemma:basicestimateun}.

%%%%%%%%%%%%%%%%%%%%%%%%%%%%%%%%%%%%%%%%%%%%%%%%%%%%%%%%
%%%%%%%%%%%%%%%%%%%%%%%%%%%%%%%%%%%%%%%%%%%%%%%%%%%%%%%%

%%%%%%%%%%%%%%%%%%%%%%%%%%%%%%%%%%%%%%%%%%%%%%%%%%%%%%%%
%%%%%%%%%%%%%%%%%%%%%%%%%%%%%%%%%%%%%%%%%%%%%%%%%%%%%%%%

\subsection{Control of the scalar products and proof of Proposition \ref{propuniq}}

%%%%%%%%%%%%%%%%%%%%%%%%%%%%%%%%%%%%%%%%%%%%%%%%%%%%%%%%
%%%%%%%%%%%%%%%%%%%%%%%%%%%%%%%%%%%%%%%%%%%%%%%%%%%%%%%%

It now remains to control the possible growth of the scalar product terms in \fref{controlnscal}. We claim:

\begin{Lemma}[A priori control of the null space]
\label{scalarproductterms}
There holds for $t$ close enough to $0$:
\be
\label{estscalarprouctmain}
\mbox{Scal}(t)\lesssim |t|^{\frac{1}{2}}|t|^2N(t).
\ee
\end{Lemma}

Let us assume Lemma \ref{scalarproductterms} and conclude the proof of Proposition \ref{propuniq}.\\

{\bf Proof of Proposition \ref{propuniq}}\\

From \fref{controlnscal}, \fref{estscalarprouctmain} and the law $\lambda_c\sim |t|$, we have for $t$ close enough to $0$:
$$N(t)\lesssim |t|^{\frac{1}{2}}N(t)+\int_t^0  \frac{N(\tau)}{\sqrt{|\tau|}}d\tau\lesssim|t|^{\frac 1 2}N(t)$$ and hence $N(t)=0$ for $t$ small enough. From the definition \fref{un70} of $N$, this yields $u=u_c$  and concludes the proof of Proposition \ref{propuniq}.\\

{\bf Proof of Lemma \ref{scalarproductterms}}\\

{\bf step 1} Approximate equation in conformal variables to the order $O(\l_c^5)$.\\

Let $v,\uv$ be defined by:
\be
\label{ps3}
u(t,x)=\frac{1}{\lambda_c(t)}v\left(t,\frac{x}{\lambda_c(t)}\right)e^{i\gamma_c(t)}, \ \ \uv(s,y)=v(s,y)e^{i\frac{b_c|y|^2}{4}}
\ee
then $\uv$ satisfies the equation:
\bee
  i\partial_s \uv+\Delta \uv-\uv+((b_c)_s+b_c^2)\frac{|y|^2\uv}{4}+\uv|\uv|^{2}=  i\left(\frac{(\l_c)_s}{\l_c}+b_c\right) \left(\Lambda \uv -ib_c\frac{|y|^2}{4}\uv\right)+(\gamma_c)_s\uv.
\eee
Starting with $u_c$, we also define $v_c$ and $\uv_c(s,y)=v_c(s,y)e^{i\frac{b_c|y|^2}{4}}$. We let $u=u_c+\ttu$ and define: 
\be
\label{defetilde}
v=v_c+\e, \ \ \uv=\uv_c+\ue, \ \ \mbox{i.e.} \ \ \ue=\e e^{i\frac{b_c|y|^2}{4}}.
\ee
Since $u_c$ satisfies \eqref{nls}, $\ue$ satisfies:
\bea
\label{ps7}
\nonumber & &  i\partial_s\ue+\Delta \ue-\ue+((b_c)_s+b_c^2)\frac{|y|^2\ue}{4} +(\uv|\uv|^{2}-\uv_c|\uv_c|^2)\\
&=&(\gamma_c)_s\ue+\ds i\left(\frac{(\l_c)_s}{\l_c}+b_c\right) \left(\Lambda \ue -ib_c\frac{|y|^2}{4}\ue\right).
\eea
\eqref{boundpmc}, \eqref{boundmodc} and \eqref{ps7} yield:
\be\label{pscl8} 
 i\partial_s\ue+\Delta \ue-\ue +(\uv|\uv|^{2}-\uv_c|\uv_c|^2)=O\left(\l_c^5(1+|y|^2)\ue+\l_c^5(1+|y|)\n \ue\right).
\ee
Let $\ue=\ue_1+i\ue_2$. We define:
$$\ue_z(s,y)=\e_z(s,y)e^{i\frac{b_c|y|^2}{4}},\,\ue_c(s,y)=\e_c(s,y)e^{i\frac{b_c|y|^2}{4}}.$$
We now expand the nonlinear term in \eqref{pscl8} as well as $\uv_c=Q+\ue_z+\ue_c$ to derive the equation at order $O(\l_c^5)$:
\be
\label{eqioheo}
-i\pa_s\ue+L(\ue)=-\psi
\ee
where $L$ is given by:
\be
\label{defnmsnkfo}
L(\ue)=L_+(\ue_1)+iL_-(\ue_2),
\ee
and where the remainder $\psi$ satisfies:
\be\label{pscl13}
\psi=O\left(\l_c^5(1+|y|^2)\ue+\l_c^5(1+|y|)\n \ue+\ue_z\ue+\ue_c\ue+\ue_z^2\ue+\ue_c^2\ue+\uv_c\ue^2+\ue^3\right).
\ee

\vspace{0.2cm}

{\bf step 2} Control of $Scal(t)$.\\

Let $f(y)=O(e^{-c|y|})$ be a smooth well localized function, then \fref{eqioheo} yields:
\be
\label{feiuyieyfei}
\frac{d}{ds}\left\{\Im(\ue,\overline{f})\right\}=-\Re(\ue,\overline{L(f)})+O((\psi,f))
\ee
with 
\bea
\label{cofeongorgh}
\nonumber |(\psi,f)| &\lesssim &\l_c^5\|\ue\|_{L^2}+\|\ue_ze^{-c|y|}\|_{L^2}\|\ue\|_{L^2}+\|\ue_c\|_{L^2}\|\ue\|_{L^2} + \|\ue_ze^{-c|y|}\|^2_{L^4}\|\ue\|_{L^2}\\
\nonumber &+ & \|\ue_c\|^2_{L^4}\|\ue\|_{L^2}+\|\uv_c\|_{L^4}\|\ue\|_{L^4}\|\ue\|_{L^2}+ \|\ue\|^2_{L^4}\|\ue\|_{L^2}\\
&\lesssim & \l_c^5\|\ue\|_{L^2}.
\eea
The control of $Scal(t)$ is now a consequence of \eqref{feiuyieyfei}, \eqref{cofeongorgh} and the structure of the null space \eqref{structurekernelpair}. 

Using \eqref{cofeongorgh}, the fact that $L_-(Q)=0$, and \eqref{feiuyieyfei} with $f=iQ$, we have:
\be\label{lundi}
\frac{d}{ds}\left\{(\ue_1,Q)\right\}=O( \l_c^5\|\ue\|_{L^2}).
\ee
We integrate this in time using the zero boundary condition from \fref{uniq} and the definition \eqref{un70} for $N(t)$:
\be
\label{lundi1}
|(\ue_1,Q)|\lesssim \int_s^{+\infty}\l_c^5\|\ue\|_{L^2}d\sigma=\int_t^{0}\l_c^3\|\ue\|_{L^2}d\tau\lesssim |t|^{\frac{7}{2}}|t|\sqrt{N(t)}.
\ee
Using \eqref{cofeongorgh}, the fact that $L_+(\Lambda Q)=-2Q$, and \eqref{feiuyieyfei} with $f=\Lambda Q$, we have:
$$\frac{d}{ds}\left\{(\ue_2,\Lambda Q)\right\}=2(\ue_1,Q)+O(\l_c^5\|\ue\|_{L^2}),$$ which time integration using \fref{lundi1} yields:
\bea
\label{lundi2}    
\nonumber |(\ue_2,\Lambda Q)| & \lesssim & \int_s^{+\infty}\left(|(\ue_1,Q)|+\l^5\|\ue\|_{L^2}\right)d\sigma\lesssim \int_t^0|\tau|^{\frac52}\sqrt{N(\tau)}d\tau\\
& \lesssim &  |t|^{\frac{5}{2}}|t|\sqrt{N(t)}.
\eea
Using \eqref{cofeongorgh}, the fact that $L_-(|y|^2Q)=-4\Lambda Q$, and \eqref{feiuyieyfei} with $f=i|y|^2Q$, we have:
$$-\frac{d}{ds}\left\{(\ue_1,|y|^2Q)\right\}=4(\ue_2,\Lambda Q)+O(\l_c^5\|\ue\|_{L^2}),$$
which together with \eqref{lundi2}  yields:
\be
\label{lundi3}
|(\ue_1,|y|^2Q)|\lesssim \int_t^0|\tau|^{\frac32}\sqrt{N(\tau)}d\tau\lesssim  |t|^{\frac{3}{2}}|t|\sqrt{N(t)}.
\ee
Using \eqref{cofeongorgh}, the fact that $L_+(\rho)=|y|^2Q$, and \eqref{feiuyieyfei} with $f=\rho$, we have:
$$\frac{d}{ds}\left\{(\ue_2,\rho)\right\}=-(\ue_1,|y|^2Q)+O(\l_c^5\|\ue\|_{L^2}),$$
which together with \eqref{lundi3} yields:
\be
\label{lundi4}
|(\ue_2,\rho)|\lesssim  \int_t^0|\tau|^{\frac12}\sqrt{N(\tau)}d\tau\lesssim  |t|^{\frac{1}{2}}|t|\sqrt{N(t)}.
\ee
Finally, \eqref{lundi1}-\eqref{lundi4} together with the fact that $\ue=\e e^{i\frac{b_c|y|^2}{4}}$ and $e^{-ib_c\frac{|y|^2}{4}}=1+O(|t||y|^2)$ imply:
$$\ds |(\e_1,Q)|+|(\e_2,\Lambda Q)|+|(\e_1,|y|^2Q)|+|(\e_2,\rho)|\lesssim |t|^{\frac{1}{2}}|t|\sqrt{N(t)}
$$
and \fref{estscalarprouctmain} is proved.\\

%%%%%%%%%%%%%%%%%

%%%%%%%%%%%%%%%%%%%%%%%%%%%%%%%%%
%%%%%%%%%%%%%%%%%%%%%%%%%%%%%%%%%


\begin{thebibliography}{10}

\bibitem{Be} Beceanu, M., A centre-stable manifold for the focussing cubic NLS in $\Bbb R^{1+3}$, Comm. Math. Phys. 280 (2008), no. 1, 145-205.

%
\bibitem{BW} Bourgain, J.; Wang, W., Construction of blowup solutions for the nonlinear Schr\"odinger equation with critical nonlinearity, Ann. Scuola Norm. Sup. Pisa Cl. Sci. (4) \textbf{25} (1997), no. 1-2, 197-215.

\bibitem{Cazenave} Cazenave, Th.; Semilinear Schr\"odinger equations, Courant Lecture Notes in Mathematics, 10, NYU, CIMS, AMS 2003.
 % 

\bibitem{CW} Cazenave, T.; Weissler, F.B., Some remarks on the nonlinear Schr{\"o}dinger equation in the critical case, Nonlinear semigroups, partial differential equations and attractors (Washington, DC, 1987), 18-29, Lecture Notes in Math., 1394, Springer, Berlin, 1989.

\bibitem{Do1}  Dodson, B., Global well posedness and scattering for the defocusing $L^2$ critical nonlinear Schr\"odinger equation when $d=2$, arxiv:1006.1375

\bibitem{Do2}  Dodson, B., Global well posedness and scattering for the defocusing $L^2$ critical nonlinear Schr\"odinger equation when $d=1$, arxiv:1010.0040

\bibitem{DR} Duyckaerts, T.; Roudenko, S., Threshold solutions for the focusing 3D cubic Schr{\"o}dinger equation. Rev. Mat. Iberoam. 26 (2010), no. 1, 1-56,

\bibitem{DM1} Duyckaerts, T.; Merle, F., Dynamic of threshold solutions for energy-critical NLS. Geom. Funct. Anal. 18 (2009), no. 6, 1787-1840.

\bibitem{DM2} Duyckaerts, T.; Merle, F., Dynamics of threshold solutions for energy-critical wave equation. Int. Math. Res. Pap. IMRP 2008, Art ID rpn002, 67 pp. 

%
\bibitem{GNN} Gidas, B.; Ni, W.M.; Nirenberg, L.,
Symmetry and related properties via the maximum principle,
Comm. Math. Phys. {\bf 68} (1979), 209-243.

%
\bibitem{GV} Ginibre, J.; Velo, G., On a class of nonlinear Schr\"odinger equations. I. The Cauchy problem, general case, J. Funct. Anal. 32 (1979), no. 1, 1-32. 

\bibitem{HR} Hillairet, M.; Rapha\"el, P., Smooth type II blow up solutions to the energy critical wave equation in dimension four,  arXiv:1010.1768.

\bibitem{Tao1} Killip, R.; Tao,T.; Visan, M.,  The cubic nonlinear
Schr{\"o}dinger equation in two dimensions with radial data.  J. Eur. Math. Soc.  11  (2009),  no. 6, 1203-1258.

\bibitem{visanetal} Killip, R.; Li, D.; Visan, M.; Zhang, X., Characterization of minimal-mass blowup solutions to the focusing mass-critical NLS, SIAM J. Math. Anal. 41 (2009), no. 1, 219-236.

\bibitem{KS2} Krieger, J.; Schlag, W., Stable manifolds for all monic supercritical focusing nonlinear Schr{\"o}dinger equations in one dimension, J. Amer. Math. Soc. 19 (2006), no. 4, 815-920.

\bibitem{KS3} Krieger, J.; Schlag, W., On the focusing critical semi-linear wave equation, Amer. J. Math. 129 (2007), no. 3, 843-913. 

\bibitem{KS1} Krieger, J.; Schlag, W. Non-generic blow-up solutions for the critical focusing NLS in 1-D, J. Eur. Math. Soc. (JEMS) 11 (2009), no. 1, 1-125.

\bibitem{KST1}  Krieger, J.; Schlag, W.; Tataru, D., Renormalization and blow up for the critical Yang-Mills problem. Adv. Math. 221 (2009), no. 5, 1445-1521. 

\bibitem{KST2} Krieger, J.; Schlag, W.; Tataru, D., Slow blow-up solutions for the $H^1(\Bbb R^3)$ critical focusing semilinear wave equation, Duke Math. J. 147 (2009), no. 1, 1-53. 

\bibitem{KST3} Krieger, J.; Schlag, W.; Tataru, D.,Renormalization and blow up for charge one equivariant critical wave maps. Invent. Math. 171 (2008), no. 3, 543-615.
 
%
\bibitem{KW} Kwong, M. K., Uniqueness of positive solutions of $\Delta u-u+u\sp p=0$ in ${R}\sp n$. Arch. Rational Mech. Anal. 105 (1989), no. 3, 243-266.%

\bibitem{LiZh} Li, D.; Zhang, X., On the rigidity of solitary waves for the focusing mass-critical NLS in dimensions $d\geq 2$. Preprint, arXiv:0902.0802.

\bibitem{martel} Martel, Y., Asymptotic $N$-soliton-like solutions of the subcritical and critical generalized Korteweg-de Vries equations, Amer. J. Math. 127 (2005), no. 5, 1103-1140.

\bibitem{MM} Matano, H.;  Merle, F.,  Classification of type I and type II behaviors for a
supercritical nonlinear heat equation.  J. Funct. Anal.  256  (2009),  no. 4,
992-1064.

\bibitem{Mc} Merle, F., Construction of solutions with exactly $k$ blow-up points for the Schr\"odinger equation with critical nonlinearity, Comm. Math. Phys. 129 (1990), no. 2, 223-240.

%
\bibitem{Mcpam} Merle, F., On uniqueness and continuation properties after blow-up time of self-similar solutions of nonlinear Schr\"odinger equation with critical exponent and critical mass, Comm. Pure Appl. Math. 45 (1992), no. 2, 203-254.

%
\bibitem{M1} Merle, F., Determination of blow-up solutions with minimal mass for nonlinear Schr\"odinger equations with critical power. Duke Math. J. 69 (1993), no. 2, 427-454.


%
\bibitem{MR1} Merle, F.; Rapha\"el, P., Blow up dynamic and upper bound on the blow up rate for critical nonlinear Schr\"odinger equation, Ann. Math. 161 (2005), no. 1, 157-222.

%
\bibitem{MR2} Merle, F.; Rapha\"el, P., Sharp upper bound on the blow up rate for critical nonlinear Schr\"odinger equation, Geom. Funct. Anal. 13 (2003), 591-642.

%
\bibitem{MR3} Merle, F.; Rapha\"el, P., On universality of blow up profile for $L^2$ critical nonlinear Schr\"odinger equation, Invent. Math. 156, 565-672 (2004).

%
\bibitem{MR4} Merle, F.; Rapha\"el, P., Sharp lower bound on the blow up rate for critical nonlinear Schr\"odinger equation, J. Amer. Math. Soc. 19 (2006), no. 1, 37-90.

%
\bibitem{MR5} Merle, F.; Rapha\"el, P.,  Profiles and quantization of the blow up mass for critical nonlinear Schr\"odinger equation, Comm. Math. Phys.  253  (2005),  no. 3, 675-704.

%
\bibitem{MR6} Merle, F.; Rapha\"el, P., On one blow up point solutions to the critical nonlinear Schr\"odinger equation, J. Hyperbolic Differ. Equ., 2 (2005), 919-962.

%
\bibitem{MR7} Merle, F.; Rapha\"el, Pierre, Blow up of the critical norm for some radial $L\sp 2$ super critical nonlinear Schr\"odinger equations, Amer. J. Math. 130 (2008), no. 4, 945-978.

\bibitem{NS1}Nakanishi, K.; Schlag, W., Global dynamics above the ground state energy for the focusing nonlinear Klein-Gordon equation, arxiv 1005.4894

\bibitem{NS2} Nakanishi, K.; Schlag, W., Global dynamics above the ground state energy for the cubic NLS equation in 3D, arxiv 1007.4025

\bibitem{PLSS} Landman, M. J.; Papanicolaou, G. C.; Sulem, C.; Sulem, P.-L., Rate of blowup for solutions of the nonlinear Schr\"odinger equation at critical dimension. Phys. Rev. A (3) 38 (1988),
no. 8, 3837-3843. 

%
\bibitem{P} Perelman, G., On the blow up phenomenon for the critical nonlinear Schr\"odinger equation in 1D, Ann. Henri. Poincaré, 2 (2001), 605-673.

\bibitem{Planchon} Planchon, F., Dispersive estimates and the 2D cubic NLS equation, J. Anal. Math. 86 (2002), 319-334.

\bibitem{PR} Planchon, F.; Rapha\"el, P., Existence and stability of the log-log blow-up dynamics for the $L^ 2$-critical nonlinear Schr\"odinger equation in a domain, Ann. Henri Poincar\'e 8 (2007), no. 6, 1177-1219.

%
\bibitem{R1} Rapha\"el, P., Stability of the log-log bound for blow up solutions to the critical nonlinear Schr\"odinger equation, Math. Ann. 331 (2005), 577-609.

%
\bibitem{R2} Rapha\"el, P., Existence and stability of a solution blowing up on a sphere for a $L^2$ supercritical nonlinear Schr\"odinger equation, Duke Math. J. 134 (2006), no. 2, 199-258.

%
\bibitem{RS} Rapha\"el, P., Szeftel, J., Existence and uniqueness of minimal mass blow up solutions to an inhomogeneous mass critical NLS, to appear in Jour. Amer. Math. Soc.

%
\bibitem{RRS} Rapha\"el, P.; Rodnianski, I., Stable blow up dynamics for the critical Wave Maps and Yang-Mills, in preparation.

\bibitem{S1} Schlag, W., Stable manifolds for an orbitally unstable nonlinear Schr{\"o}dinger equation. Ann. of Math. (2) 169 (2009), no. 1, 139-227.

\bibitem{Str} Strichartz, R. S., Restrictions of Fourier transforms to quadratic surfaces and decay of solutions of wave equations, Duke Math. J. 44 (1977), no. 3, 705-714.

%
\bibitem{SS} Sulem, C.; Sulem, P.L., The nonlinear Schr\"odinger equation. Self-focusing and wave collapse. Applied Mathematical Sciences, 139. Springer-Verlag, New York,
1999.

\bibitem{W1}  Weinstein, M.I., Nonlinear Schr\"odinger equations and sharp interpolation estimates, Comm. Math. Phys.
{\bf 87} (1983), 567-576.

\bibitem{W2}  Weinstein, M.I., Modulational stability of ground 
states of nonlinear Schr\"odinger equations, SIAM J. Math. Anal.
{\bf 16} (1985), 472-491.

%
\bibitem{ZS} Zakharov, V.E.; Shabat, A.B., Exact theory
of two-dimensional self-focusing and one-dimensional self-modulation
of waves in non-linear media, Sov. Phys. JETP 34 (1972),
62-69.

\end{thebibliography}
\end{document}